\documentclass{article}
\usepackage[a4paper, margin=2.5cm]{geometry}
\usepackage{lmodern}
\usepackage{algorithm}
\usepackage{comment}
\usepackage{algpseudocode}
\usepackage{amsthm}
\usepackage{graphicx}
\usepackage{placeins}
\usepackage{float}
\usepackage{amsmath}
\usepackage[nocompress]{cite}
\usepackage{amssymb}
\usepackage[title]{appendix}
\usepackage[
    colorlinks=true,
    citecolor=blue,
    linkcolor=black,
    urlcolor=blue
]{hyperref}

\theoremstyle{remark}

\theoremstyle{plain}
\newtheorem{theorem}{Theorem}[section]
\newtheorem{lemma}[theorem]{Lemma}
\newtheorem{proposition}[theorem]{Proposition}

\theoremstyle{definition}
\newtheorem{assumption}[theorem]{Assumption}
\newcommand{\NErrorEval}{500}

\numberwithin{equation}{section}

\numberwithin{algorithm}{section}

\title{On multi-fidelity methods for a tumor growth model with uncertainties}

\author{
Yu Feng$^{2}$,
Liu Liu$^{1,*}$,
Te Qi$^{1}$,
and Huimin Yu$^{1}$
}

\date{}

\begin{document}

\maketitle

\begin{center}
{\small
$^{1}$ The Chinese University of Hong Kong, Hong Kong, P. R. China\\
$^{2}$ The Great Bay University, Dongguan 523808, P. R. China
}
\end{center}

\renewcommand{\thefootnote}{}
\footnotetext{
\raggedright
$^*$ Corresponding author: Liu Liu (liuliu@cuhk.edu.hk).

E-mail addresses:
1155275626@link.cuhk.edu.hk (Huimin Yu),
fengyu@gbu.edu.cn (Yu Feng),
1155145997@link.cuhk.edu.hk (Te Qi).

This work was supported by the National Key R\&D Program of China
(2021YFA1001200), the Ministry of Science and Technology of China,
the General Research Fund of the Research Grants Council of Hong Kong
(14303022, 14301423, 14307125, 14310726), and the National Natural Science
Foundation of China (NSFC Youth Program, 12501669).
}
\begin{abstract}
In this work, we develop a hierarchical multi-fidelity method for uncertainty quantification of porous-medium equation (PME) models of tumor growth with moving interfaces. With an accurate asymptotic-preserving scheme for the PME model as our high-fidelity solver, we choose an inexpensive level-set method for the limiting Hele--Shaw equation as its low-fidelity solver. In contrast to existing approaches, in this work high-fidelity parameter samples are selected by the pivoted Cholesky followed by a residual-based farthest-point sampling step to better resolve interface variations of the solution. Moreover, an empirical error estimator is given to assess the resulting bi- and tri-fidelity reconstructions. Several numerical experiments under various uncertainty settings have demonstrated an accurate and efficient recovery of tumor morphology and uncertainty statistics.
\end{abstract}

\textbf{Keywords:} porous medium equation, uncertainty quantification, multi-fidelity methods, residual-based farthest-point sampling

\noindent\textbf{MSC:}
65M75, 65C20, 35R35, 92C50.

\section{Introduction}
Tumor growth models couple nutrient diffusion, cell proliferation, and tissue mechanics. Classical tumor-growth models \cite{chen2003free,friedmanreitich2001symmetry,friedman2007mathematical,greenspan1972models} treat the tumor as an incompressible fluid flowing through a porous medium, governed by equations of Hele--Shaw type. In recent years, connections with porous-medium-type equations (PME) have substantially advanced the mathematical understanding of these models. In the PME models, pressure and density are related by the constitutive law \(p=\frac{m}{m-1}\rho^{m-1}\), where the exponent \(m\) controls the stiffness of the pressure law. The stiff-pressure limit \(m\to\infty\) connects porous-medium tumor growth models to a Hele--Shaw free-boundary description \cite{perthame2014hele,perthame2014traveling}. This correspondence has motivated further analysis of incompressible limits and interface dynamics \cite{david2021free,Feng2024jde,jacobs2023tumor,kim2018porous,liu2025three,tong2025convergence}. It also underlies numerical methods for the density and free-boundary models \cite{Feng2024unified,liu2018accurate,liu2021toward}.

Uncertainty quantification (UQ) is an important component of tumor growth modeling, since biological heterogeneity, parameter variability, and limited observability can affect model predictions \cite{pareschi2021introduction}. Standard numerical approaches to UQ include Monte Carlo sampling, stochastic collocation methods \cite{babuvska2007stochastic,nobile2008sparse}, and generalized polynomial chaos--stochastic Galerkin (gPC-SG) methods \cite{babuska2004galerkin}. Applications to tissue and tumor growth modeling include Bayesian inversion \cite{falco2023quantifying,Feng2024unified} and stochastic asymptotic-preserving schemes \cite{jiang2026efficient}. For the PME model, increasing stochastic dimension and the stiffening pressure law as \(m\to\infty\) make accurate resolution of perturbation-induced interface variations computationally demanding.

Multi-fidelity (MF) methods reduce the cost of repeated high-fidelity (HF) simulations by combining a limited number of HF evaluations with inexpensive low-fidelity (LF) models \cite{peherstorfer2016optimal,peherstorfer2018survey}. For solution reconstruction, approximation coefficients computed in the LF solution space are transferred to a reduced HF basis \cite{narayan2014mf,zhu2014siamjuq}. In practice, LF models often arise from coarse discretizations or asymptotic reductions, such as the Hele--Shaw limit of the PME. Multi-fidelity methods have also been developed for kinetic and related models \cite{DimarcoLiuPareschiZhu2025,dimarco2019multi,lin2025sisc,LiuPareschiZhu2022jcp,LiuZhu2020jcp}. 

A crucial step in multi-fidelity methods is selecting representative samples from the LF solution space. Pivoted Cholesky decomposition and rank-revealing QR factorization identify dominant approximation directions in this space \cite{harbrecht2012cholesky,Hampton2018jcp,narayan2014mf,zhu2014siamjuq}. Accurate recovery of tumor morphology and interface evolution requires capturing both the overall shape and fine-scale variations near the moving boundary \cite{gil2001convergence,kim2018porous,vazquez2007porous,feng2023tumor}. To balance global representativeness with sensitivity to local interface variations, residual-based farthest-point sampling (RFPS) combines fitting residuals with the separation criterion of farthest-point sampling (FPS) \cite{eldar1997farthest,gonzalez1985clustering}.

In this work, we develop a multi-fidelity framework for
porous-medium-type tumor growth models with multi-dimensional
uncertain parameters. The high-fidelity PME model is solved by an asymptotic-preserving scheme that is stable and accurate in both moderate and large-\(m\) regimes
\cite{liu2018accurate,liu2021toward}. The low-fidelity models are
constructed from coarse-grid PME discretizations and a level-set method for the Hele--Shaw limit
\cite{alam2025thresholding,osher1988fronts,sethian1999level}. We construct both bi-fidelity (BF) and tri-fidelity
(TRI) approximations by using the hierarchy among fine-grid PME,
coarse-grid PME, and Hele--Shaw models.

The main contributions and novelty of this work are summarized as follows:
\begin{itemize}
    \item Within the multi-fidelity framework, we adopt the inexpensive Hele--Shaw
    limit as our low-fidelity solver and make modifications of existing bi-fidelity method for kinetic models in \cite{LiuZhu2020jcp,LinLiu2026mms} through a two-stage sample-selection procedure to better capture the interface-related
    features of the PME model with uncertain parameters.
    \item We analyze our new enrichment stage in the
    important samples selection procedure and develop an empirical error estimator to assess the reliability of multi-fidelity reconstructions by using a limited
    number of high-fidelity simulations.
\end{itemize}

The rest of the paper is organized as follows. Section~2 introduces the
porous-medium tumor growth model and its asymptotic Hele--Shaw limit.
Section~3 presents the multi-fidelity reconstruction framework, the
two-stage sample-selection procedure, and the fidelity solvers used in
the numerical tests. Section~4 analyzes the stability
and error of the enrichment procedure and develops an empirical error
estimator based on the error structure.
Section~5 reports numerical experiments under several
uncertainty scenarios, including perturbations in initial conditions,
nutrient-related parameters, and interface shapes. Section~6 concludes
the paper.

\section{Tumor growth model}
We consider mechanically driven tumor growth models of porous-medium type
and their incompressible Hele--Shaw limit in two spatial dimensions,
$x\in\mathbb{R}^2$. These models describe tumor expansion driven by pressure-induced motion, while proliferation is
regulated by nutrient availability
\cite{liu2019analysis,perthame2014hele,perthame2014traveling}.
Uncertainty is introduced through a random parameter
$z\in I_z\subset\mathbb{R}^{d_z}$ endowed with probability density
$\pi(z)\ge0$. For each fixed $z$, the system is deterministic.

Let $\rho=\rho(x,t,z)$ denote the tumor-cell density and $c=c(x,t,z)$ the nutrient concentration. Let
\[
D^+(t,z)=\{x\in\mathbb{R}^2:\rho(x,t,z)>0\}
\]
denote the tumor region. The governing system couples mass conservation
with Darcy's law:
\begin{equation}\label{eq:main_clean}
\partial_t\rho+\nabla\cdot(\rho u)=g(x,t,\rho,z),
\qquad u=-\nabla p.
\end{equation}
The pressure is $p=p(\rho)$, and $m$ is a parameter in the constitutive law below. The proliferation
term is given by
\[
g(x,t,\rho,z)=G_0 c(x,t,z)\rho,
\]
where $G_0=G_0(z)>0$ is the proliferation coefficient; its dependence on $z$ is suppressed in the notation.

The quasi-steady in-vivo nutrient model \cite{liu2019analysis}
balances diffusion with tumor-cell consumption $\lambda(z)\rho c$
and vascular supply outside the tumor. The consumption coefficient
$\lambda(z)$ and background concentration $c_B(z)$ are nonnegative;
the exterior supply coefficient is normalized to one. The nutrient field satisfies
\begin{equation}\label{eq:c_vivo}
\begin{cases}
-\Delta c + \lambda(z)\rho c = 0,
& x \in D^+(t,z),\\[4pt]
-\Delta c + c = c_B(z),
& x \in \mathbb{R}^2\setminus D^+(t,z).
\end{cases}
\end{equation}
The nutrient concentration and its normal derivative are continuous
across the tumor boundary, and $c\to c_B(z)$ as $|x|\to\infty$.

\vspace{2mm}
\noindent\textbf{Porous medium model. }
To close the system, we assume the constitutive relation between pressure
and density is given by the power law
\[
p=\frac{m}{m-1}\rho^{m-1},
\qquad m\ge2.
\]
Let $f(x,z)$ be the prescribed initial tumor-cell density, with $0\le f\le1$.
Then \eqref{eq:main_clean} reduces to the porous-medium equation
\begin{equation}\label{Pm_clean}
\begin{cases}
\partial_t\rho-\dfrac{m}{m-1}\nabla\cdot(\rho\nabla\rho^{m-1})
=G_0c(x,t,z)\rho,\\[4pt]
\rho(x,0,z)=f(x,z).
\end{cases}
\end{equation}
This equation is a degenerate diffusion--reaction system in flux form; see \cite{vazquez2007porous} for the porous-medium diffusion theory.
As $m$ increases, the pressure law becomes increasingly stiff, giving
rise to the congestion constraint $\rho_\infty\le1$ in the limit.
\vspace{2mm}

\noindent\textbf{Hele--Shaw limit. }
The formal stiff-pressure limit of \eqref{Pm_clean} as $m\to\infty$
is a Hele--Shaw free-boundary model. The limit density and pressure
$(\rho_\infty,p_\infty)$ satisfy the Hele--Shaw complementarity graph
\begin{equation}
\label{eqn:HSgraph}
0\le\rho_\infty\le1,\qquad p_\infty\ge0,\qquad
p_\infty(1-\rho_\infty)=0.
\end{equation}
Passing formally to the limit in the mass balance \eqref{Pm_clean} gives
\begin{equation}\label{P_limit}
\begin{cases}
\partial_t\rho_{\infty}
-\nabla\cdot(\rho_{\infty}\nabla p_{\infty})
=
G_0 c(x,t,z)\rho_{\infty},\\[4pt]
\rho_{\infty}(x,0,z)=f(x,z),
\end{cases}
\end{equation}
in the weak sense, where $c$ now denotes the limiting nutrient field.
For smooth finite-$m$ solutions, the pressure law gives
$\partial_t p=(m-1)p(\Delta p+G_0c)+|\nabla p|^2$.
Dividing by $m-1$ and formally passing to the limit yields the
complementarity relation
\begin{equation}
\label{eqn:comple-condi}
p_{\infty}(\Delta p_\infty+G_0c(x,t,z))=0.
\end{equation}
Rigorous versions of this limit and the complementarity relation
are established for mechanical growth models
\cite[Sec.~2]{perthame2014hele} and related nutrient-coupled models
\cite{david2021free}, including the quasi-steady setting in
\cite[Appendix~A, Theorem~A.2]{Feng2024jde}.

For initial data of the form $f(x,z)=\chi_{D_0(z)}$, where $\chi_A$ denotes the characteristic function of the set $A$, the limiting model admits patch-type solutions. In this regime, the limiting density remains binary-valued and evolves as $\rho_\infty(x,t,z)=\chi_{D(t,z)}$, where
$D(t,z)\subset\mathbb{R}^2$ denotes the tumor domain at time $t$.
The limiting mass balance \eqref{P_limit}, together with
\eqref{eqn:HSgraph}, gives
\begin{equation}\label{eq:HS_closure}
\begin{cases}
-\Delta p_\infty=G_0c(x,t,z),
& x\in D(t,z),\\[4pt]
p_\infty=0,
& x\in\partial D(t,z),
\end{cases}
\end{equation}
The interface normal velocity is
$V=-\nabla p_\infty\cdot n$ on $\partial D(t,z)$, where $n$ is the
unit outward normal. In the patch regime the nutrient model
\eqref{eq:c_vivo} reduces to
\begin{equation}
\begin{cases}
-\Delta c + \lambda(z) c = 0,
& x \in D(t,z),\\[4pt]
-\Delta c + c = c_B(z),
& x \in \mathbb{R}^2\setminus D(t,z).
\end{cases}
\end{equation}

The finite-$m$ PME resolves diffuse density dynamics, while its
Hele--Shaw limit describes a saturated region with a moving interface.
These complementary models provide the basis of the fidelity hierarchy
used below.

\section{Multi-fidelity approximation framework}
\label{sec:mf_framework}

We seek low-cost reconstructions of the uncertain density field
using paired model solutions and coefficient transfer
\cite{narayan2014mf,zhu2014siamjuq}. The finite-$m$ PME provides the
HF reference, and its Hele--Shaw limit provides a level-set (LS)
interface surrogate. A coarse-grid PME (CMPME) supplies an additional lower-cost
model: it preserves the finite-$m$ equation and parameter dependence
while reducing the spatial and temporal discretization costs
\cite[Sec.~3.1.2]{Hampton2018jcp}. Solver details are given in
Section~\ref{sec:solver}.

For each parameter $z\in I_z$, let $U^H(z),U^L(z),U^M(z)
\in\mathbb R^{N_h}$ denote the discretized HF, low-fidelity (LF), and
intermediate-fidelity density fields. In BF, $U^L$ is supplied by either
CMPME or LS; in TRI, $U^L$ is supplied by CMPME and $U^M$ by LS.
All fields refer to the same output time and are represented on a common
Cartesian grid of $N_h=N_XN_Y$ points, with $N_X$ and $N_Y$ points in the
two coordinate directions. CMPME outputs are interpolated onto this grid
where needed. We use $\langle u,v\rangle=u^Tv$ and denote the BF and TRI
reconstructions by $U^B$ and $U^T$.

In the \emph{offline stage}, we evaluate the lower-cost model at
training parameters $\Gamma_N=\{z_1,\ldots,z_N\}\subset I_z$ and
select representative parameters using Section~\ref{subsec:Two-stage greedy strategy}.
Each parameter in the union of the training and enrichment pools has
a unique label. The ordered list $\gamma=(\gamma(1),\ldots,\gamma(K))$
records the selected labels, and $z_{\gamma(j)}$ denotes the parameter
with label $\gamma(j)$. This order is shared by the paired basis matrices.
Here $K=|\gamma|$ is the number of paired basis fields. At the selected parameters we compute the HF solutions and form
$\mathcal U^Y(\gamma)=\operatorname{span}\{U^Y(z_{\gamma(j)})\}_{j=1}^{K}$
and the corresponding HF space $\mathcal U^H(\gamma)$, with $Y=L$ for BF
or $Y=M$ for TRI.

In the \emph{online stage}, for a new parameter $z$ we solve the
chosen lower-cost model, fit $U^Y(z)$ in the selected basis
of $\mathcal U^Y(\gamma)$ to obtain coefficients $c_j^Y(z)$, and use
the same coefficients to combine the paired HF snapshots. Thus the
selected parameters and HF solutions are reused for new inputs;
Algorithm~\ref{alg:adaptive_bi_fidelity} gives the reconstruction;
$\operatorname{cond}$ and $\operatorname{tr}$ denote the spectral condition
number and matrix trace. Coefficient transfer relies on compatible parameter dependence across
fidelities \cite{narayan2014mf,zhu2014siamjuq}; the corresponding
error condition and estimator are developed in
Section~\ref{sec:stability_empirical}.

\begin{algorithm}[!ht]
\caption{Bi- and tri-fidelity reconstruction}
\label{alg:adaptive_bi_fidelity}
\begin{algorithmic}[1]
\State \textbf{Offline stage:}
\State Draw $N$ training parameters $\Gamma_N\subset I_z$ according to $\pi$.
\State Obtain the selected parameter labels $\gamma$ from the
two-stage selection in Section~\ref{subsec:Two-stage greedy strategy}.
\State Construct $\mathcal U^Y(\gamma)$ from the selected lower-fidelity
solutions, with $Y=L$ for BF and $Y=M$ for TRI.
\State Run the HF model at the selected parameters, set $K=|\gamma|$, and
construct $\mathcal U^H(\gamma)$. A reconstruction requires $K\ge1$ and a nonzero lower-cost basis.
\State \textbf{Online stage:}
\State For a new parameter $z\in I_z$, compute the lower-fidelity solution $U^Y(z)$.
\State Form
$(G^Y)_{\ell j}=\langle U^Y(z_{\gamma(\ell)}),U^Y(z_{\gamma(j)})\rangle$
and $g_j^Y=\langle U^Y(z),U^Y(z_{\gamma(j)})\rangle$.
\State If $K>1$ and $\operatorname{cond}(G^Y)>10^{12}$,
replace $G^Y$ by $G^Y+10^{-8}\operatorname{tr}(G^Y)I/K$.
Solve $G^Y\mathbf c^Y(z)=\mathbf g^Y(z)$ for $\mathbf c^Y(z)\in\mathbb R^K$.
\State Transfer the same coefficients to the HF solutions:
\[
U^F(z)=\sum_{j=1}^{K} c_j^Y(z)U^H(z_{\gamma(j)}),
\qquad (F,Y)=(B,L)\text{ or }(T,M).
\]
\State \Return $U^F(z)$.
\end{algorithmic}
\end{algorithm}

\subsection{Two-stage sample selection}
\label{subsec:Two-stage greedy strategy}

Let $\Gamma_N=\{z_i\}_{i=1}^N$ be the training set and
$S_L=[f_1,\ldots,f_N]$, where $f_i=U^L(z_i)$.
Stage~I selects up to $K_0$ parameters by pivoted Cholesky
\cite{harbrecht2012cholesky,zhu2014siamjuq}.
The selected LF solutions form the initial basis $\Phi_0$ for BF;
for TRI, the selected CMPME parameters are used to construct
$\Phi_0$ from LS solutions.

Stage~II selects additional parameters from the candidate pool
$\mathcal C=\{\zeta_i\}_{i=1}^{N_C}$.
If the pool consists of the remaining training parameters, then
$N_C=N-K_0$; otherwise, $N_C$ denotes the size of the prescribed
enrichment pool. Let
$F_{\mathcal C}=[g_1,\ldots,g_{N_C}]$, where
$g_i=U^Y(\zeta_i)$, with $Y=L$ for BF and $Y=M$ for TRI.

The enrichment criterion combines spatially weighted fitting
residuals with Euclidean separation between candidate fields.
For the current basis $\Phi$ and fixed $\delta>0$, define
\[
C_\delta=(\Phi^T\Phi+\delta I)^{-1}\Phi^TF_{\mathcal C},
\qquad
\mathcal R=F_{\mathcal C}-\Phi C_\delta.
\]
Using fixed nonnegative spatial weights $w_\ell$, with
$W=\operatorname{diag}(w_1,\ldots,w_{N_h})$ and
$S_w=\sum_{\ell=1}^{N_h}w_\ell>0$, set
\[
r_i^{\rm sel}
=\left(
\frac{\sum_{\ell=1}^{N_h}w_\ell\mathcal R_{\ell i}^2}
{S_w}
\right)^{1/2},
\qquad
J(\Phi)=
\left[\frac{1}{N_C}\sum_{i=1}^{N_C}
(r_i^{\rm sel})^2\right]^{1/2}.
\]
The coefficient fit and separation distance use the ordinary
inner product, while $\omega$ balances the normalized residual
and distance scores.

Starting from $J_0=J(\Phi_0)$, enrichment stops when
$J(\Phi)\le\theta J_0$ ($0<\theta<1$), after at most
$K_{\rm add}$ accepted additions, or upon exhaustion of the
candidate pool. Only accepted additions enter the HF basis construction. Algorithm~\ref{alg:two_stage_greedy} summarizes
the procedure, with screening and update details given in
Appendix~\ref{app:sampling}.

\begin{algorithm}[!ht]
\small
\caption{Two-stage sample selection by pivoted Cholesky and RFPS}
\label{alg:two_stage_greedy}
\begin{algorithmic}[1]
\Require $S_L$, nonempty $F_{\mathcal C}$, fixed spatial weights $W$, $\delta>0$;
$K_0\ge1$, $K_{\rm add}$, $\omega\in[0,1]$, $\theta\in(0,1)$
\Ensure Selected labels $\gamma$, basis $\Phi$, and stopping status
\Statex \textbf{Stage I: pivoted-Cholesky initialization}
\State Apply pivoted Cholesky to $S_L^TS_L$ until $K_0$ selections or
numerical rank exhaustion; obtain $\gamma$, set $K_0\gets|\gamma|$.
\State Form the nonempty initial basis $\Phi=\Phi_0$ from the corresponding
$U^Y$ fields ($Y=L$ for BF; $Y=M$ for TRI).
\Statex \textbf{Stage II: residual-based farthest-point sampling (RFPS)}
\State Find $j_a\in\arg\min_j\|g_j-\Phi(:,a)\|_2$, $a=1,\dots,K_0$.
\State Initialize $d_i=\min_a\|g_i-g_{j_a}\|_2$,
$\mathcal E=\varnothing$, $k=0$, $J_0=J(\Phi)$.
\While{true}
  \State Compute $C_\delta$, $\mathcal R$, $r_i^{\rm sel}$ and $J(\Phi)$.
  \State \textbf{if} $J(\Phi)\le\theta J_0$, \textbf{break} (target reached).
  \State \textbf{if} $k=K_{\rm add}$ or $|\mathcal E|=N_C$,
  \textbf{break} (budget or pool exhausted).
  \State Normalize residuals to $[0,1]$ and distances by their maximum
  over the fixed pool (Appendix~\ref{app:sampling}).
  \State \textbf{RFPS selection:} $i^\star\in\arg\max_{i\notin\mathcal E}
  \{\omega\widetilde r_i+(1-\omega)\widetilde d_i\}$.
  \State $\mathcal E\gets\mathcal E\cup\{i^\star\}$.
  \State \textbf{if} the residual or correlation checks fail, \textbf{continue}.
  \State Append $g_{i^\star}$ to $\Phi$ and its parameter label to $\gamma$;
  $k\gets k+1$.
  \State \textbf{FPS update:} $d_i\gets\min\{d_i,\|g_i-g_{i^\star}\|_2\}$
  for $i\notin\mathcal E$.
\EndWhile
\State \Return $\gamma,\Phi$ and stopping status; $K_0+k$ locations are selected.
\end{algorithmic}
\end{algorithm}

\subsection{High- and lower-fidelity solvers}
\label{sec:solver}

The fine-grid finite-$m$ PME is solved by an asymptotic-preserving
(AP) method \cite{liu2018accurate,liu2021toward}. CMPME uses its coarse-grid
counterpart, and the level-set (LS) solver approximates the Hele--Shaw
limit. Either CMPME or LS supplies the lower-cost model for BF. For TRI
(Test~2(b)), we use the CMPME--LS--PME hierarchy: CMPME is used for
initial sample selection, while LS provides complementary interface
information for RFPS enrichment and coefficient recovery.

\subsubsection{An asymptotic-preserving method for the high-fidelity model}

Our HF solver follows the prediction--correction AP framework of
\cite{liu2018accurate,liu2021toward}, which remains consistent with the
Hele--Shaw limit as $m\to\infty$.
\vspace{1mm}

\textbf{Velocity--density formulation. }
Introducing the velocity
$u=-\nabla p(\rho)=-\frac{m}{m-1}\nabla\rho^{m-1}$, the density equation
can be written as
$\partial_t\rho+\nabla\!\cdot(\rho u)=G_0c(x,t,z)\rho$. The AP solver
evolves the coupled density--velocity equations and enforces the
constitutive relation
$u=-\frac{m}{m-1}\nabla\rho^{m-1}$ through a correction step. This gives
the constrained system
\begin{equation}\label{eq:HF_relaxed_short}
\begin{cases}
\partial_t \rho + \nabla\!\cdot(\rho u)
=G_0c(x,t,z)\rho,\\[3pt]
\partial_t u
=
m\nabla\!\left(
\rho^{m-2}
\big(\nabla\!\cdot(\rho u)-G_0c(x,t,z)\rho\big)
\right),\\[3pt]
u=-\dfrac{m}{m-1}\nabla\rho^{m-1}.
\end{cases}
\end{equation}
\vspace{1mm}

\textbf{Prediction--correction splitting. }
Given $(\rho^n,u^n)$, the prediction step advances the first two
equations of \eqref{eq:HF_relaxed_short} to obtain $(\rho^\ast,u^\ast)$.
The correction step keeps
$\rho^{n+1}=\rho^\ast$ and restores the consistency relation
\[
u^{n+1}
=
-\frac{m}{m-1}\nabla(\rho^{n+1})^{m-1}.
\]
This splitting separates the transport--reaction update from the
nonlinear diffusion constraint. The discrete update formulas are given
in Appendix~\ref{app:PME}.

\subsubsection{A level-set solver for the Hele--Shaw low-fidelity model}

We approximate the Hele--Shaw interface by a first-order
Hamilton--Jacobi level-set method \cite{osher1988fronts,sethian1999level},
using the level-set discretization discussed in
\cite[Sec.~3]{alam2025thresholding}.

Let $\phi(x,t,z)$ be a level-set function such that
$D(t,z)=\{x:\phi(x,t,z)<0\}$ and
$\partial D(t,z)=\{x:\phi(x,t,z)=0\}$. The interface evolves according to
the Hamilton--Jacobi equation
\begin{equation}\label{eq:level_set}
\phi_t+V(x,t,z)|\nabla\phi|=0,
\end{equation}
where $V$ is an extended normal velocity field. At each time step, given the current interface
$D^k=\{x:\phi^k(x)<0\}$, we solve the Hele--Shaw pressure equation
\begin{equation}\label{eq:HS_pressure_discrete}
-\Delta p_\infty^k=G_0c(x,t^k,z),\quad x\in D^k,
\qquad
p_\infty^k=0,\quad x\in\partial D^k.
\end{equation}
The interface normal velocity is computed from Darcy's law as
$V_n^k=-\nabla p_\infty^k\cdot n$, where
$n=\nabla\phi^k/|\nabla\phi^k|$ is the outward unit normal.

The normal speed is extended off the interface by solving
$\nabla\phi^k\cdot\nabla V^k=0$ along its normals, as in
\cite[Sec.~3.3]{alam2025thresholding}. We then advance the level-set
equation explicitly with a first-order upwind Godunov discretization
and periodically reinitialize $\phi$ near the interface using fast sweeping
\cite[Secs.~3.4--3.5]{alam2025thresholding}. For coefficient transfer,
the LS field is represented by the density indicator
$\rho_{\mathrm{LS}}(x,t,z)=\chi_{\{\phi(x,t,z)<0\}}$ on the common
density grid.

\section{Stability and error assessment}
\label{sec:stability_empirical}

We analyze the regularized fitting and enrichment in
Algorithm~\ref{alg:two_stage_greedy}, followed by the error of the HF
reconstruction. On the common grid of $N_h=N_XN_Y$ points, use the scaled
spatial norm $\|U\|_{s,h}=\|U\|_2/N_h$ for $U\in\mathbb R^{N_h}$, as in
the numerical error measure in Section~\ref{sec:numerical}. For vectors,
$\|\cdot\|_2$ is the Euclidean norm; for matrices, $\|\cdot\|_2$ and
$\|\cdot\|_F$ denote the spectral and Frobenius norms. The RMS aggregation
of the scaled errors bounds their arithmetic mean on the same sample set.

\subsection{Fitting stability and enrichment analysis}
\label{subsec:regularized_analysis}

Recall the selection metric $W=\operatorname{diag}(w_1,\ldots,w_{N_h})$
and $S_w=\sum_\ell w_\ell>0$ from Section~\ref{subsec:Two-stage greedy strategy}.
Let $w_{\max}=\max_\ell w_\ell>0$. For a fixed selection basis
$\Phi\in\mathbb R^{N_h\times r}$ and fixed $\delta>0$, set
$A_\delta=(\Phi^T\Phi+\delta I)^{-1}\Phi^T$,
$C_\delta=A_\delta F_{\mathcal C}$, and
$Q_\delta=I_{N_h}-\Phi A_\delta$, so that
$\mathcal R=Q_\delta F_{\mathcal C}$.

\begin{lemma}[Fixed-basis stability]
\label{lem:ridge_stability}
For $\delta>0$, the regularized fit is unique. With $\Phi$ fixed,
a perturbation $\Delta F$ of $F_{\mathcal C}$ satisfies
\begin{equation}
\label{eq:ridge_stability}
\|\Delta C_\delta\|_F
\le \|A_\delta\|_2\|\Delta F\|_F,
\qquad
\|A_\delta\|_2
=\max_j\frac{\sigma_j}{\sigma_j^2+\delta}
\le\frac{1}{2\sqrt\delta}.
\end{equation}
The $\sigma_j$ are the singular values of $\Phi$.
For a fixed HF basis $\Psi$ whose columns are paired with those of
$\Phi$ in the same order, the reconstruction perturbation is at most
$\|\Psi\|_2\|A_\delta\|_2\|\Delta F\|_F$ in Frobenius norm.
\end{lemma}

The proof uses the singular-value representation of the Tikhonov-regularized
inverse \cite{hansen2010discrete}. For a fixed retained basis, the same
estimate applies to the final coefficient solve with its actual positive
regularization parameter. If that solve is unregularized, full column rank
gives the gain $1/\sigma_{\min}(\Phi)$.

\begin{proposition}[Quantitative decrease after accepted enrichment]
\label{prop:accepted_descent}
For fixed $\delta>0$ and $F_{\mathcal C}$, define
$\mathcal V_\delta(\Phi)=\min_C\{
\|F_{\mathcal C}-\Phi C\|_F^2+\delta\|C\|_F^2\}$.
Let $a=g_i$ be an accepted column with pre-enrichment residual
$r_i^{\rm sel}>0$, and set $\Phi_+=[\Phi,a]$ and
$\alpha_a=S_w(r_i^{\rm sel})^2/w_{\max}$. Then
\begin{equation}
\label{eq:objective_decrease}
\mathcal V_\delta(\Phi)-\mathcal V_\delta(\Phi_+)
=\frac{\|F_{\mathcal C}^TQ_\delta a\|_2^2}
{\delta+a^TQ_\delta a}
\ge\frac{\alpha_a^2}{\delta+\alpha_a}>0.
\end{equation}
For every selection basis,
$J(\Phi)^2\le w_{\max}\mathcal V_\delta(\Phi)/(N_CS_w)$.
\end{proposition}

The descent identity follows from the Sherman--Morrison rank-one
inverse formula \cite[Eq.~(2)]{hager1989updating}. The proofs of
Lemma~\ref{lem:ridge_stability} and Proposition~\ref{prop:accepted_descent}
are given in Appendix~\ref{app:analysis_proofs}. The accepted addition
decreases the upper envelope of $J^2$; Algorithm~\ref{alg:two_stage_greedy}
checks the weighted residual itself at each stopping test.

\subsection{Conditional reconstruction error control}
\label{subsec:conditional_error}

At a reconstruction prefix, let $r$ be the number of paired basis
columns, $\Phi_r^Y,\Phi_r^H$ the lower-cost and HF bases, and $C_r$
the coefficient matrix actually used in Algorithm~\ref{alg:adaptive_bi_fidelity}.
The coefficient-recovery model is $Y=L$ for BF and $Y=M$ for TRI.
Here $C_r\in\mathbb R^{r\times N_C}$ contains the coefficients for
all enrichment candidates, and
$H_{\mathcal C}=[U^H(\zeta_1),\ldots,U^H(\zeta_{N_C})]$
contains their paired HF fields in the same order. Define
\begin{equation}
\label{eq:error_rms_analysis}
E_r=\frac{1}{N_C}\sum_{i=1}^{N_C}
\|H_{\mathcal C}(:,i)-\Phi_r^HC_r(:,i)\|_{s,h},
\qquad
E_{2,r}=\frac{\|H_{\mathcal C}-\Phi_r^HC_r\|_F}{N_h\sqrt{N_C}}.
\end{equation}
The subscript $2$ denotes RMS aggregation, not a different spatial
normalization. Put
$\mathsf Y_r=[F_{\mathcal C},\Phi_r^Y]$ and
$\mathsf H_r=[H_{\mathcal C},\Phi_r^H]$.

\begin{assumption}[Weighted Gram compatibility]
\label{ass:paired_compatibility}
For specified constants $\beta_r,\varepsilon_r\ge0$, the paired
dictionaries satisfy
\begin{equation}
\label{eq:gram_compatibility}
G_{H,r}:=\frac{\mathsf H_r^T\mathsf H_r}{N_h^2}
\preceq
\beta_r G_{Y,r}+\varepsilon_r I,
\qquad
G_{Y,r}:=\frac{\mathsf Y_r^TW\mathsf Y_r}{S_w}.
\end{equation}
\end{assumption}

Here $W$ and $S_w$ are defined in Section~\ref{subsec:Two-stage greedy strategy},
$I=I_{N_C+r}$, and $A\preceq B$ means that $B-A$ is positive semidefinite.
This assumption adapts the Gram-defect comparison in
\cite[Eq.~(9) and Lemma~2]{Hampton2018jcp} to the weighted selection metric
and the scaled HF norm, with $\varepsilon_r$ accounting for model differences.

\begin{theorem}[Reconstruction error and residual stopping]
\label{thm:stop_error}
Under Assumption~\ref{ass:paired_compatibility}, set
$J_r=\|W^{1/2}(F_{\mathcal C}-\Phi_r^YC_r)\|_F/
\sqrt{N_CS_w}$ and $q_r=1+\|C_r\|_F^2/N_C$. Then
\begin{equation}
\label{eq:gram_error_bound}
E_r\le E_{2,r}\le
\bigl(\beta_rJ_r^2+\varepsilon_rq_r\bigr)^{1/2}.
\end{equation}
Suppose selection stops with $s$ columns and
$J(\Phi_s)\le\theta J_0$, and Algorithm~\ref{alg:adaptive_bi_fidelity}
uses the corresponding $K=s$ paired columns satisfying the Gram compatibility
condition \eqref{eq:gram_compatibility} for $r=K$. Let $C_{\delta,s}$ be the
selection coefficients and define
the change in the lower-cost fit on the same candidate pool
\begin{equation}
\label{eq:fit_mismatch}
\Delta_{\rm fit}=
\frac{\|W^{1/2}(\Phi_sC_{\delta,s}-\Phi_K^YC_K)\|_F}
{\sqrt{N_CS_w}}.
\end{equation}
For the final retained basis,
\begin{equation}
\label{eq:successful_stop_bound}
E_K\le E_{2,K}\le
\bigl[\beta_K(\theta J_0+\Delta_{\rm fit})^2+
\varepsilon_Kq_K\bigr]^{1/2}.
\end{equation}
\end{theorem}

\begin{proof}
Let $D_r=\begin{bmatrix}I_{N_C}\\-C_r\end{bmatrix}\in
\mathbb R^{(N_C+r)\times N_C}$. The HF and lower-cost residuals are
$\mathsf H_rD_r$ and $\mathsf Y_rD_r$, respectively. Applying
\eqref{eq:gram_compatibility} to the columns of $D_r$ and summing gives
\[
 \frac{\|\mathsf H_rD_r\|_F^2}{N_h^2N_C}
 \le \beta_r\frac{\|W^{1/2}\mathsf Y_rD_r\|_F^2}{S_wN_C}
     +\varepsilon_r\frac{\|D_r\|_F^2}{N_C}.
\]
Since $\|D_r\|_F^2=N_C+\|C_r\|_F^2$, this is the squared RMS bound in
\eqref{eq:gram_error_bound}; the mean bound follows from Cauchy--Schwarz.
For the final basis, add and subtract $\Phi_sC_{\delta,s}$ in the
lower-cost residual to obtain $J_K\le J(\Phi_s)+\Delta_{\rm fit}$.
The stopping target and \eqref{eq:gram_error_bound} at $r=K$ then give
\eqref{eq:successful_stop_bound}.
\end{proof}

The bound \eqref{eq:successful_stop_bound} connects the Gram-matrix
comparison to RFPS stopping, with $\Delta_{\rm fit}$ accounting for
differences between the selection and reconstruction coefficient fits.
For fixed $\beta_r$, the smallest nonnegative admissible defect is
$\varepsilon_r=\max\{0,\lambda_{\max}(G_{H,r}-\beta_rG_{Y,r})\}$,
where $\lambda_{\max}$ is the largest eigenvalue. Evaluating this defect
requires HF fields on the paired pool; the following empirical estimator
uses a small calibration set instead.

\subsection{Error estimator}
\label{subsection:stable_empirical}

Following the error decomposition in \cite[Sec.~2.2.4]{Gao2020cmame},
we construct an empirical estimator using the error outside and within the
retained HF space. Let $P_r^Y$ and $P_r^H$ be the
orthogonal projectors in the ordinary inner product onto the retained
lower-cost and HF spaces.
For a finite paired sample set $\mathcal A$, use
\eqref{eq:error_rms_analysis} with $\mathcal C$ replaced by $\mathcal A$
to define $E_r(\mathcal A)$ and $E_{2,r}(\mathcal A)$, and set
\begin{equation}
\label{eq:projection_transfer_scales}
\begin{aligned}
D_{Y,r}(\mathcal A)&=
\frac{\|(I-P_r^Y)F_{\mathcal A}\|_F}{N_h\sqrt{|\mathcal A|}},
&D_{H,r}(\mathcal A)&=
\frac{\|(I-P_r^H)H_{\mathcal A}\|_F}{N_h\sqrt{|\mathcal A|}},\\
T_{H,r}(\mathcal A)&=
\frac{\|P_r^HH_{\mathcal A}-\Phi_r^HC_{r,\mathcal A}\|_F}
{N_h\sqrt{|\mathcal A|}}.
\end{aligned}
\end{equation}
The matrices $F_{\mathcal A}$ and $H_{\mathcal A}$ contain the paired
snapshots, and $C_{r,\mathcal A}$ contains the actual reconstruction coefficients.
The residual splits as
$H_{\mathcal A}-\Phi_r^HC_{r,\mathcal A}=(I-P_r^H)H_{\mathcal A}+
(P_r^HH_{\mathcal A}-\Phi_r^HC_{r,\mathcal A})$.
The two terms lie in orthogonal spaces, giving the exact identity
\begin{equation}
\label{eq:orthogonal_error_identity}
E_{2,r}(\mathcal A)^2
=D_{H,r}(\mathcal A)^2+T_{H,r}(\mathcal A)^2.
\end{equation}
On a fixed paired calibration set $\Gamma_{\rm cal}$ with
$D_{Y,r}(\Gamma_{\rm cal})>0$, define
$\widehat a_r=D_{H,r}(\Gamma_{\rm cal})/
D_{Y,r}(\Gamma_{\rm cal})$ and
$\widehat b_r=T_{H,r}(\Gamma_{\rm cal})$.
The resulting estimate of the RMS error $E_{2,r}(\mathcal A)$ is
\begin{equation}
\label{eq:calibrated_estimator}
\widehat E_r(\mathcal A)=
\left[\widehat a_r^{\,2}D_{Y,r}(\mathcal A)^2+
\widehat b_r^{\,2}\right]^{1/2}.
\end{equation}
The factor $\widehat a_r$ estimates the HF projection-error scale
from the lower-cost residual, and $\widehat b_r$ estimates the error
within the HF space. The estimator relies on the calibration set being
representative of the target set. Section~\ref{sec:numerical} gives the
calibration and resampling settings.

\section{Numerical experiments}
\label{sec:numerical}

All tests consider tumor growth on
$\Omega=[-2.5,2.5]^2$ up to $T=1$, with a fine Cartesian grid of
$N_X=N_Y=101$ points ($\Delta x=0.05$). We impose
$\rho u\cdot n=0$ and $c=c_B(z)$ on $\partial\Omega$ for the density
flux and nutrient concentration, respectively. CMPME fields are
transferred to the HF grid by bilinear interpolation when needed;
LS indicators are evaluated directly on that grid.

CMPME uses $\Delta x_2=0.1$ with parabolic time-step scaling
$\Delta t\sim\Delta x^2$ for both PME resolutions. For $m=8$, the fine-grid
PME time step is $1.5\times10^{-3}$ and the CMPME time step is $6\times10^{-3}$.
For $m=64$, these are $5\times10^{-4}$ and $2\times10^{-3}$, respectively.
LS uses the fine grid and $\Delta t_{\mathrm{LS}}=0.03$.

The uncertain input is $z=(z_1,\ldots,z_{d_z})$, with independent
components $z_j\sim\mathcal U(-1,1)$; we set $d_z=5$ in all tests.
In all tests, the LF training set $\Gamma_N$ contains $N=500$
parameters. In all constructions, pivoted Cholesky selects the initial
$K_0=15$ HF sample locations.The enrichment pool $\mathcal C$ consists of the remaining
training parameters after the initial selection.
RFPS may add up to $K_{\rm add}=5$ further HF sample locations, giving
an HF basis budget of at most 20.

An independent Monte Carlo ensemble
$\Gamma_{N_s}=\{z^{(i)}\}_{i=1}^{N_s}$, $N_s=1200$, is used to compare
reconstructed and HF statistics at $T=1$. For any of the density fields
$U$, the pointwise sample mean and standard deviation are
\[
 \overline U=\frac1{N_s}\sum_{i=1}^{N_s}U(T,z^{(i)}),\qquad
 S_U=\left[\frac1{N_s}\sum_{i=1}^{N_s}
 (U(T,z^{(i)})-\overline U)^2\right]^{1/2},
\]
which estimate $\mathbb E[U]$ and $\mathrm{Std}[U]$, respectively.
In all tests, affine random multipliers are clipped from below at $10^{-6}$
before constructing the physical parameters, using the same preprocessing
for all fidelities and reference computations.
The PME initial densities are kept below unity for numerical stability.

All RFPS experiments use the same Stage-II acceptance checks in
Appendix~\ref{app:sampling} and the residual stopping criterion in
Algorithm~\ref{alg:two_stage_greedy}. The normalized residual and
separation scores have equal weight $\omega=0.5$, and the stopping ratio
is $\theta=0.5$. The number of accepted additions is determined by this
stopping procedure within the prescribed budget. The Gaussian weights are
$w_\ell=\exp[-(x_\ell^2+y_\ell^2)/(2h_w^2)]$, where
$(x_\ell,y_\ell)$ is the $\ell$th common-grid point and $h_w=0.15$
is the Gaussian width. The weights emphasize central density
variations, while the unweighted distance term in RFPS incorporates
separation of the full fields.

For the field reconstruction error, we use a fixed, independently sampled evaluation set
$\Gamma_{\rm eval}=\{z_i^{\rm eval}\}_{i=1}^{N_{\rm eval}}$ with
$N_{\rm eval}=\NErrorEval$. At each evaluation parameter, the HF
reference and the BF or TRI approximation use the \emph{same} $z_i^{\rm eval}$;
the latter is reconstructed from the selected HF basis using the
corresponding lower-cost solution. The reported mean scaled spatial error is
\begin{equation}
\label{eq:mean_error}
E=
\frac{1}{N_{\rm eval}}
\sum_{i=1}^{N_{\rm eval}}
\left(
\frac{1}{N_XN_Y}
\sqrt{
\sum_{j=1}^{N_X}
\sum_{k=1}^{N_Y}
\left|
U^H(T,x_j,y_k,z_i^{\rm eval})-U^F(T,x_j,y_k,z_i^{\rm eval})
\right|^2
}
\right),
\end{equation}
where $U^F$ denotes the BF or TRI reconstruction. Equivalently, $E$
is the sample mean of
$\|U^H(T,z_i^{\rm eval})-U^F(T,z_i^{\rm eval})\|_{s,h}$.

After basis selection, the estimator is calibrated at each
retained basis size $r$ using the same $|\Gamma_{\rm cal}|=20$ matched
lower-cost/HF samples. We draw 2000 paired bootstrap resamples of the
calibration indices \cite[Sec.~2]{efron1979bootstrap}, recomputing
$\widehat a_r$ and $\widehat b_r$ in each resample while keeping the
bases and evaluation residuals fixed. The calibration ratio requires
a positive lower-cost projection residual in each resample. The plotted
estimator is the empirical $0.95$ quantile of the resulting
$\widehat E_r(\Gamma_{\rm eval})$ values. This quantile estimates the RMS
error scale and is compared with the mean error $E$ in \eqref{eq:mean_error}.
Calibration HF evaluations are counted separately from the HF basis budget.

The HF Monte Carlo reference ensemble was generated using the
Supercomputing Center of Wuhan University. The reported per-realization
solver timings refer to individual forward solves and were measured
separately in MATLAB R2021a on Windows 10 Home, using an Intel Core i5-8300H
CPU (2.30 GHz, 4 cores and 8 logical processors) with 16 GB of RAM.

\subsection{Test 1: Interface perturbations with coupled parameter uncertainty}

This test compares an unweighted LF-residual greedy baseline with
RFPS enrichment after the same pivoted-Cholesky initialization, under
coupled perturbations of the initial interface and model parameters.
CMPME supplies the LF solutions, and the PME solver with $m=8$ supplies
the HF solutions.
For each realization $z_i=(z_{1,i},\ldots,z_{d_z,i})$, the perturbed
interface is defined by
\[
R_{\mathrm{petal}}(x,y;z_i)
=
R_{0,i}\bigl[1+A_i\cos(n_{\mathrm{petal}}\varphi)\bigr],
\]
where \(R_0=1\), \(A_0=0.25\), \(n_{\mathrm{petal}}=6\), and \(\varphi=\operatorname{atan2}(y,x)\) is the polar angle.

The stochastic parameters are given by
\[
\lambda_i=\lambda_0(1+z_{1,i}),\quad
c_{B,i}=c_{B,0}(1+z_{1,i}),\quad
G_{0,i}=G_0(1+z_{1,i}),
\]
\[
R_{0,i}
=R_0\Bigl(1+\sum_{k=1}^{d_z}\frac{z_{k,i}}{2k}\Bigr),
\quad
A_i
=A_0\Bigl(1+\sum_{k=1}^{d_z}\frac{z_{k,i}}{2k}\Bigr),
\]
with baseline parameters
$\lambda_0=50$, $c_{B,0}=20$, and $G_0=0.5$.
The parameters $\lambda_i$, $c_{B,i}$, and $G_{0,i}$
depend on $z_1$, while $R_{0,i}$ and $A_i$ depend on
all random components.
The initial density is defined by
\[
\rho_0(x,y;z_i)=
\begin{cases}
\rho_{00,i}, & \sqrt{x^2+y^2}<R_{\mathrm{petal}}(x,y;z_i),\\
0, & \text{otherwise},
\end{cases}
\qquad
\rho_{00,i}
=
\min\bigl(0.75(1+0.5z_{1,i}),\,0.95\bigr).
\]

\begin{figure}[!htbp] 
    \centering
    \includegraphics[width=0.4\textwidth]{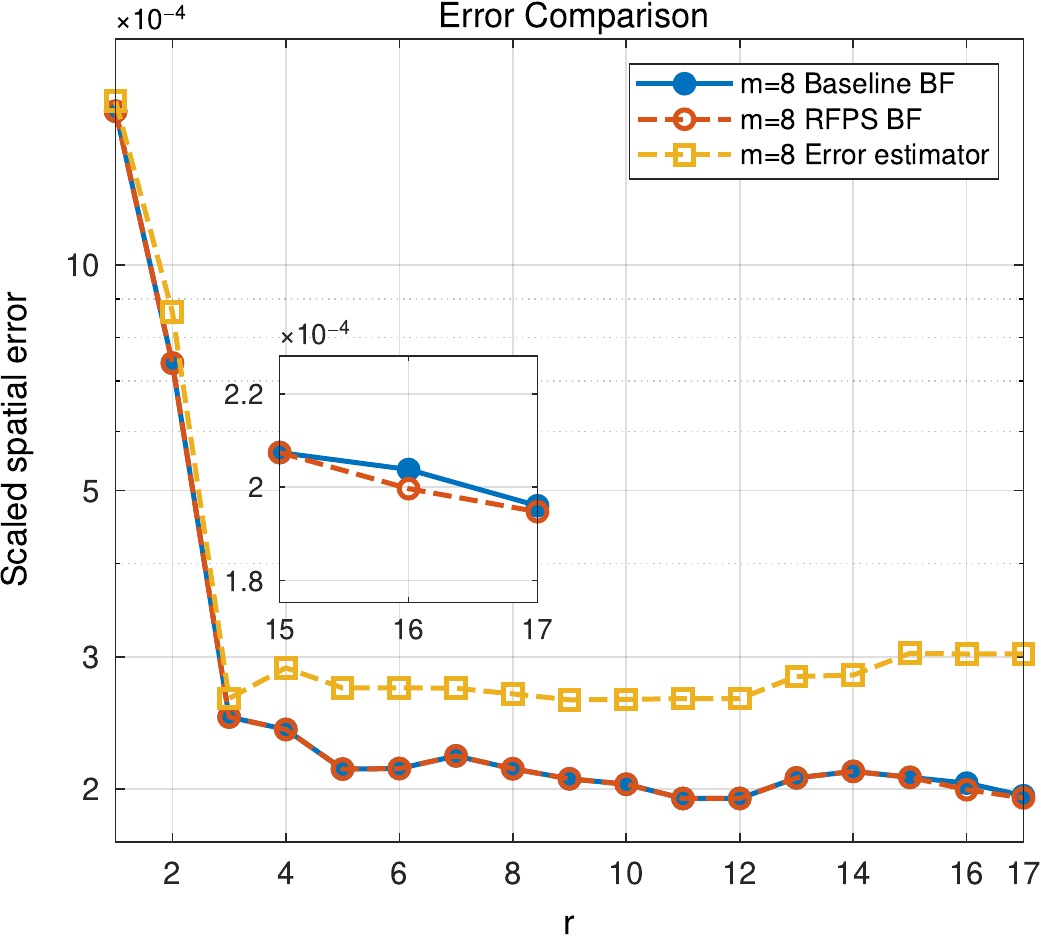}
    \caption{
Test 1: Mean scaled spatial error on the parameter-matched evaluation
set versus retained HF basis size for $m=8$. Baseline BF and RFPS-BF
are shown up to 17 fields, with the error estimator for RFPS-BF.
}
    \label{fig:error_test1}
\end{figure}

Figure~\ref{fig:error_test1} shows a rapid error reduction over the
first few basis fields, followed by smaller fluctuations. Starting from $K_0=15$, RFPS
reaches $J\le0.5J_0$ after two additions, giving $K=17$. The additions
slightly reduce the mean reconstruction error, and the empirical estimator
remains above the observed RFPS-BF error over the displayed prefixes.

\begin{figure}[!htbp] 
    \centering
    \includegraphics[width=0.9\textwidth]{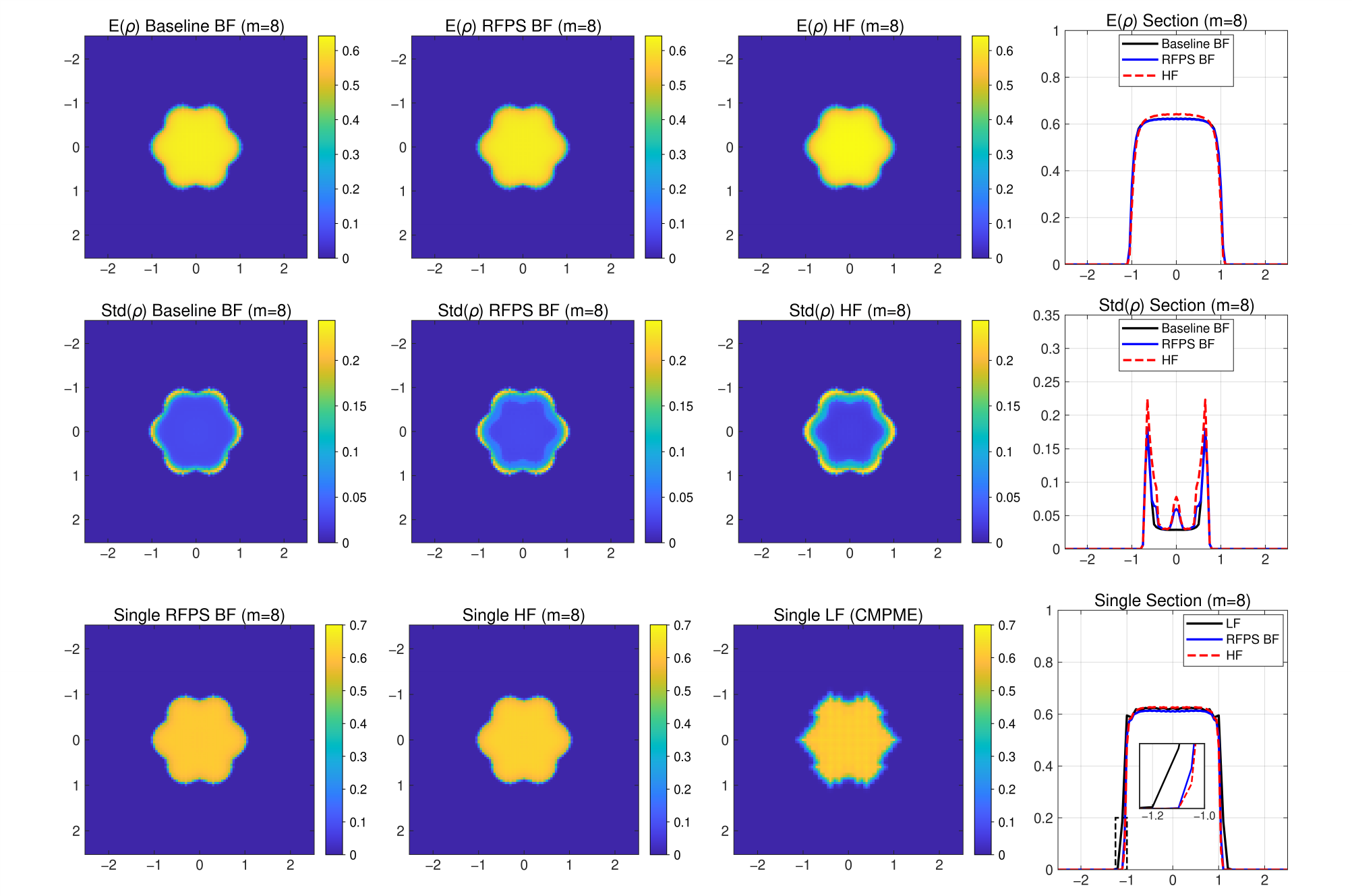}
\caption{
Test 1: Comparison of BF approximations for $m=8$, using 18 baseline
and 17 RFPS basis fields. The first row
shows the mean density field and the second row shows the standard
deviation, with columns corresponding to baseline BF, RFPS-BF, the HF
reference, and one-dimensional slices. The third row shows one
representative realization at a fixed parameter $z$, including the RFPS-BF reconstruction, the HF
reference, the CMPME solution recomputed with the same initial field
and model parameters, and the corresponding slice.
Slices are taken at \(y=0\) for the mean field and the representative
realization, and at $y=0.65$ for the standard deviation.
}
    \label{fig:test1_mean_std}
\end{figure}

Both reconstructions in Figure~\ref{fig:test1_mean_std} recover the
six-lobed mean morphology. RFPS improves the interior standard-deviation
profile while retaining the peaks associated with boundary uncertainty.
For the representative realization, CMPME captures the overall tumor
shape, but its density transition near the interface differs from HF.
The RFPS slice closely follows the HF profile in this region, as the inset shows.

For Test~1 with $m=8$, the reported CMPME and HF PME times are
approximately $8.349$ s and $318.393$ s per realization, respectively,
giving an HF-to-CMPME cost ratio of about $38.14$.
Using CMPME for online evaluations therefore reduces the cost per forward
solve relative to direct HF sampling.

\subsection{Test 2: Nutrient-related uncertainty under BF and TRI frameworks}
\label{sec:test1a}

This test examines nutrient-related uncertainty. Part~(a) compares
baseline and RFPS-enriched BF reconstructions using LS as the LF model.
Part~(b) compares CMPME-BF, LS-BF, and the RFPS-enriched TRI construction.

For each realization, the parameters are given by
\[
\lambda_i=\lambda_0\Bigl(1+\sum_{k=1}^{d_z}\frac{z_{k,i}}{2k}\Bigr),
\qquad
c_{B,i}=c_{B,0}\Bigl(1+\sum_{k=1}^{d_z}\frac{z_{k,i}}{2k}\Bigr),
\]
\[
G_{0,i}=G_0(1+z_{1,i}),
\qquad
R_{0,i}=R_0(1+z_{1,i}),
\]
with \(\lambda_0=50\), \(c_{B,0}=20\), \(G_0=0.5\), and \(R_0=0.45\).
The initial condition is radially symmetric:
\[
\rho_0(x,y;z_i)=
\begin{cases}
0.95, & \sqrt{x^2+y^2}<R_{0,i},\\
0, & \text{otherwise}.
\end{cases}
\]
We use $m=8$ and $m=64$.

\vspace{3mm}
\noindent\textbf{Test 2(a): BF comparison under two sampling strategies.}\par\nobreak
Fine-grid PME supplies the HF model. Both strategies start from
the same 15-sample initialization.
RFPS reaches $J\le0.5J_0$ after three additions, giving $K=18$.
Both strategies use 18 basis fields for the mean and standard-deviation
plots, and their error curves show the corresponding prefixes.
The baseline sequentially selects the candidate maximizing
$\|\Phi^Tg_i\|_2$ and uses $G^Y+10^{-10}I$ for coefficient recovery;
RFPS uses Algorithm~\ref{alg:adaptive_bi_fidelity} for recovery.

\begin{figure}[!htbp] 
    \centering
    \includegraphics[width=0.75\textwidth]{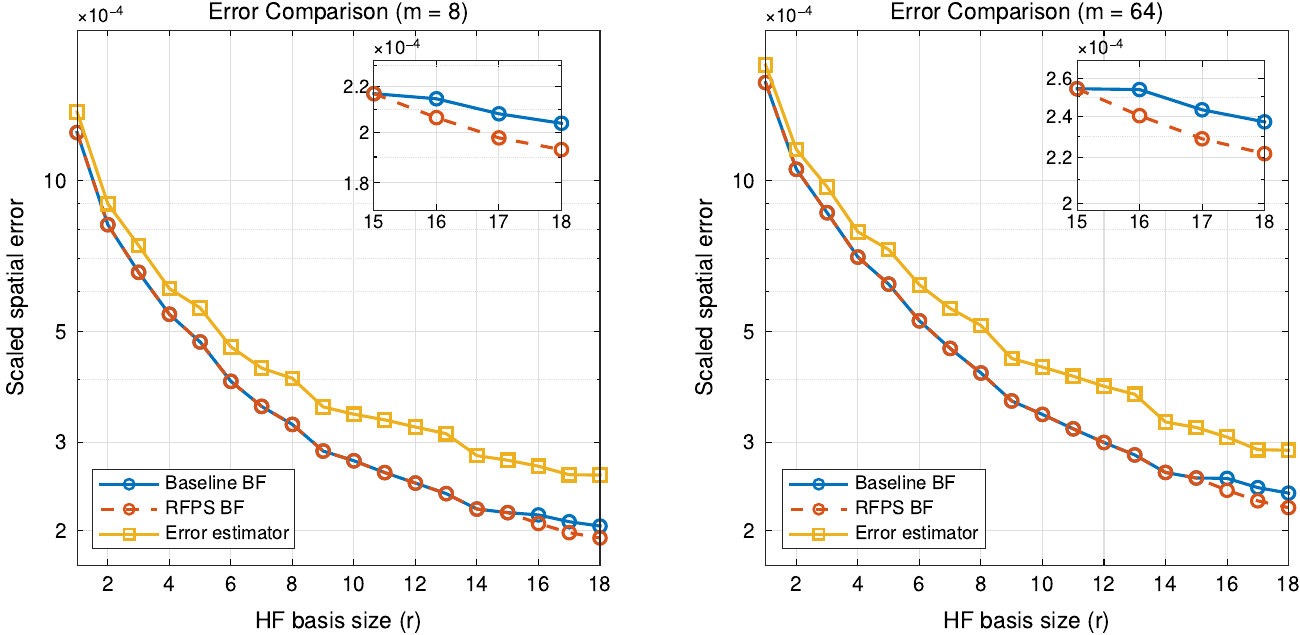}
    \caption{
Test 2(a): Mean scaled spatial error on the parameter-matched
evaluation set versus retained HF basis size. Baseline BF and RFPS-BF
for $m=8$ and $m=64$, with the RFPS-BF error estimator.
}
    \label{fig:error_test2a}
\end{figure}

\begin{figure}[!htbp] 
    \centering
    \includegraphics[width=0.82\textwidth]{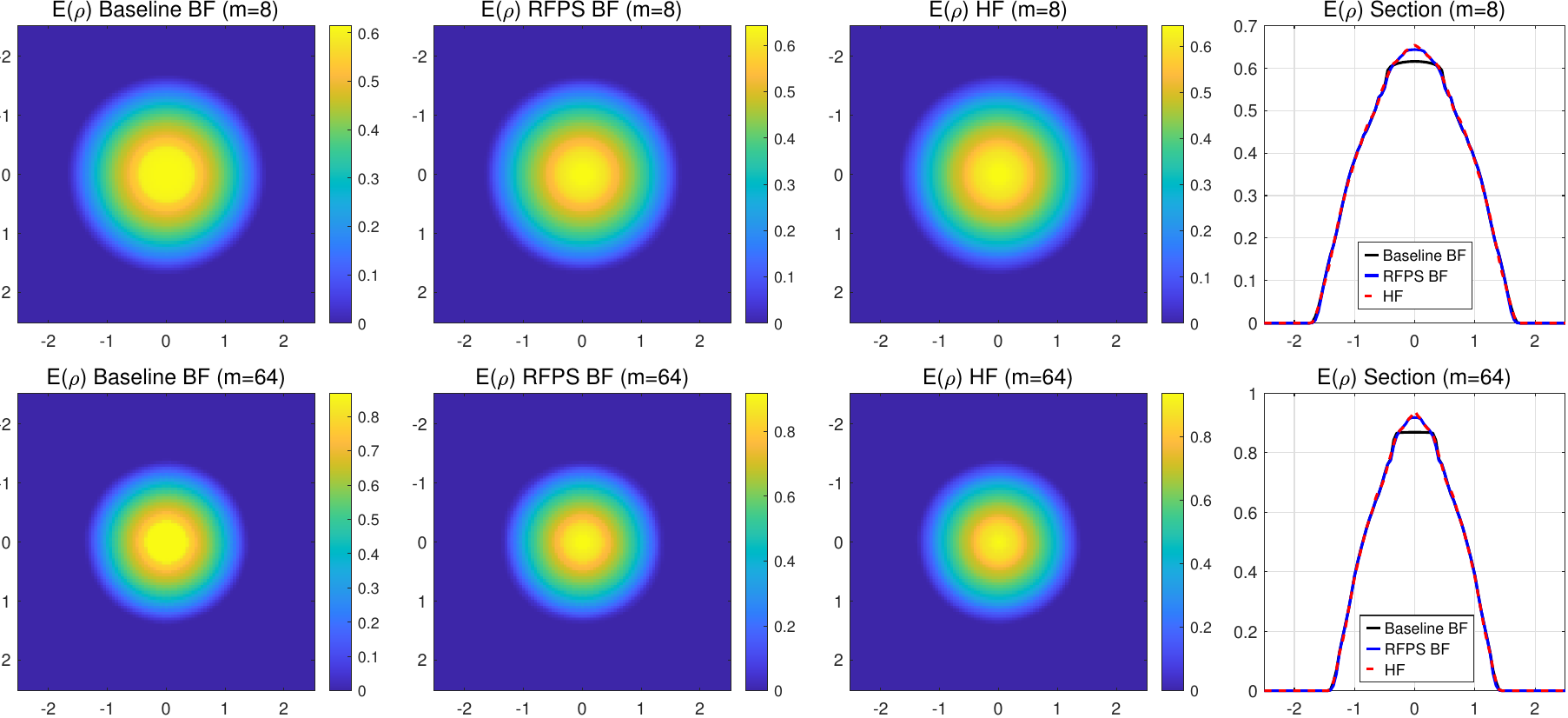}
    \caption{
Test 2(a): Mean field of \(\rho\) under BF approximation.
Top: \(m=8\); bottom: \(m=64\).
Columns show baseline BF, RFPS-BF, HF reference, and one-dimensional
slices at \(y=0\).
}
    \label{fig:test2a_mean}
\end{figure}

\begin{figure}[!htbp] 
    \centering
    \includegraphics[width=0.82\textwidth]{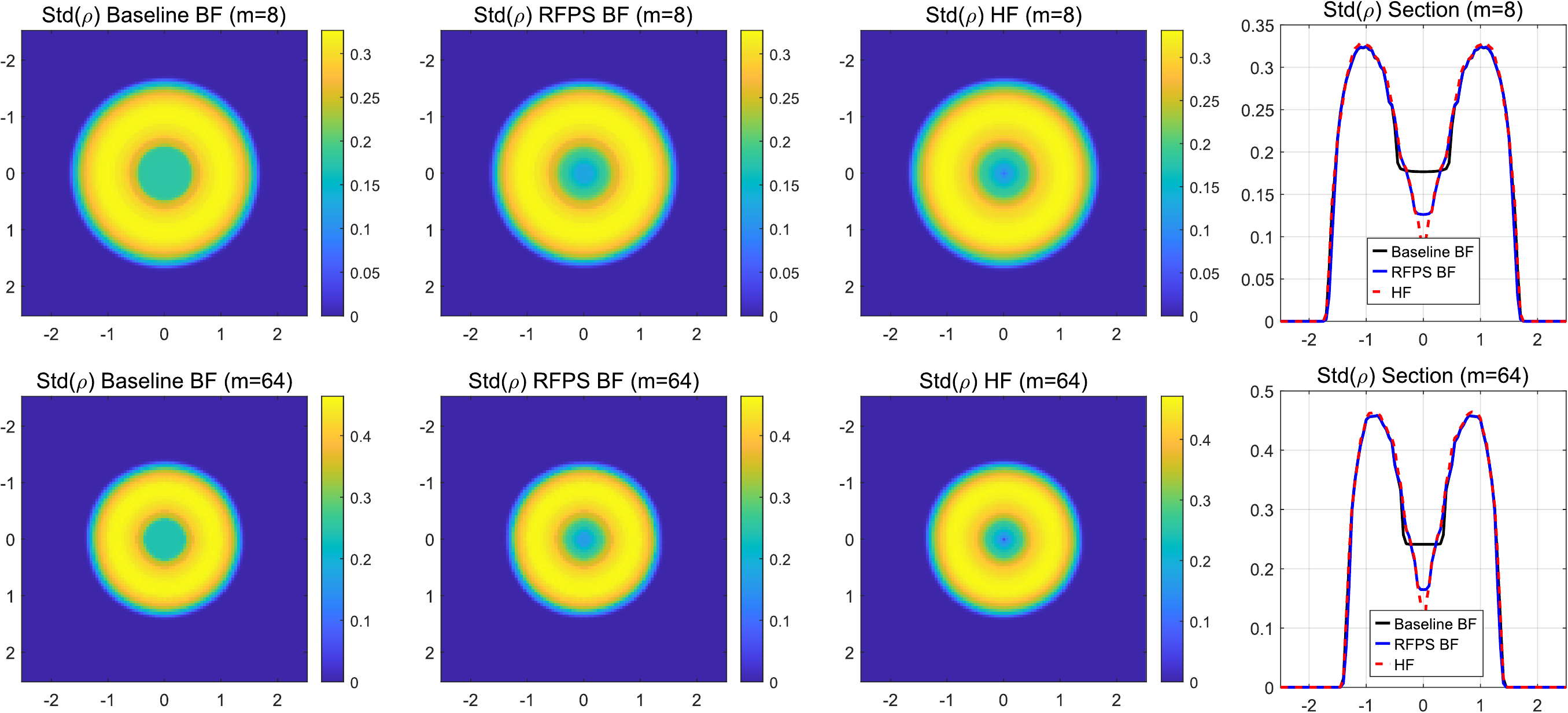}
    \caption{
Test 2(a): Standard deviation of \(\rho\) under BF approximation.
Top: \(m=8\); bottom: \(m=64\).
Columns show baseline BF, RFPS-BF, HF reference, and one-dimensional
slices at \(y=0\).
}
    \label{fig:test2a_std}
\end{figure}

Figure~\ref{fig:error_test2a} shows rapid error decay during
initialization, followed by slower decay as the basis grows. The RFPS-BF reconstruction yields a lower final error
for both $m=8$ and $m=64$, and the estimator follows the
observed RFPS-BF error from above.
Figure~\ref{fig:test2a_mean} shows that RFPS recovers the radial mean profile,
including the higher central density for $m=64$.
In Figure~\ref{fig:test2a_std}, the annular region of large standard deviation
reflects uncertainty in the tumor boundary; RFPS improves the central
depression while retaining the peripheral peaks.

\vspace{3mm}
\noindent\textbf{Test 2(b): TRI comparison under RFPS sampling.}\par\nobreak

The model parameters remain those of part~(a). We compare
LS-BF, CMPME-BF, and TRI. In TRI, CMPME selects the first 15 locations;
LS supplies the basis for RFPS enrichment and coefficient recovery.
TRI reaches $J\le0.5J_0$ after three additions, giving $K=18$.
The two enriched BF references also use 15 initial and three additional
fields, so their error curves are compared at matching HF basis sizes.

\begin{figure}[!htbp] 
    \centering
    \includegraphics[width=0.75\textwidth]{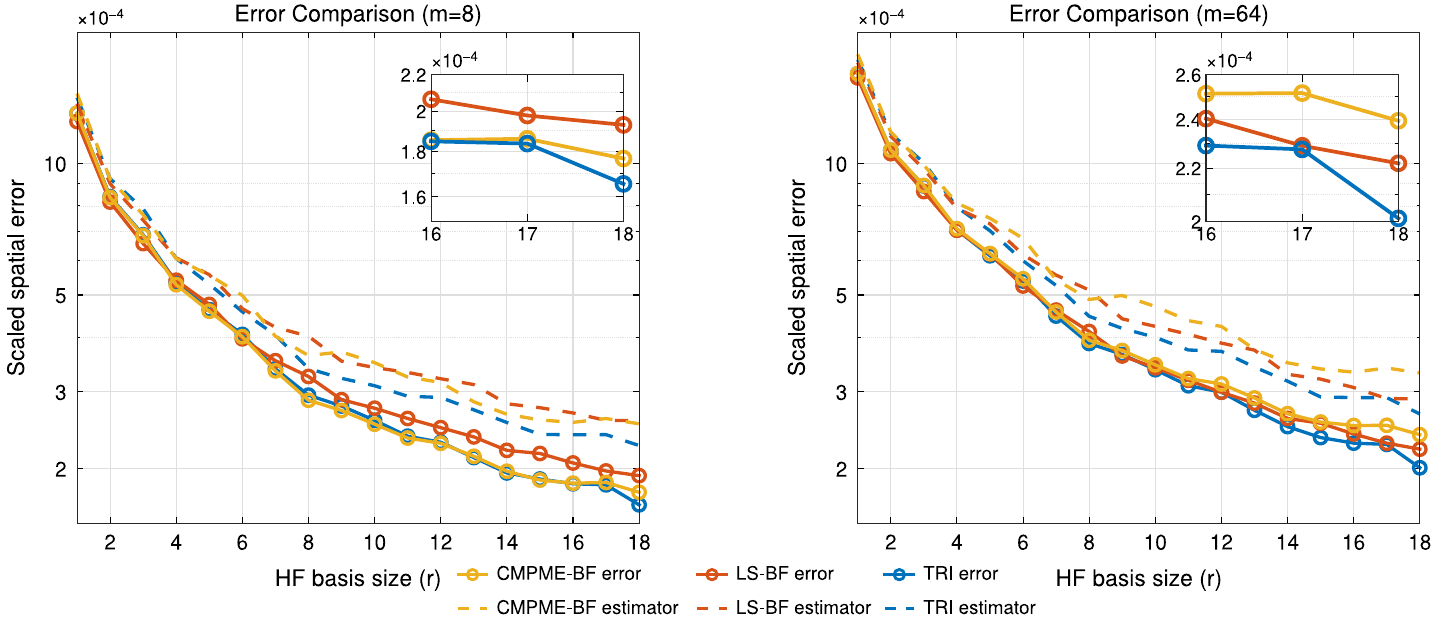}
    \caption{
Test 2(b): Mean scaled spatial error on the parameter-matched
evaluation set versus retained HF basis size. CMPME-BF, LS-BF, and TRI
for $m=8$ and $m=64$, with their corresponding error estimators.
}
    \label{fig:error_test2ab}
\end{figure}

\begin{figure}[!htbp] 
    \centering
    \includegraphics[width=0.82\textwidth]{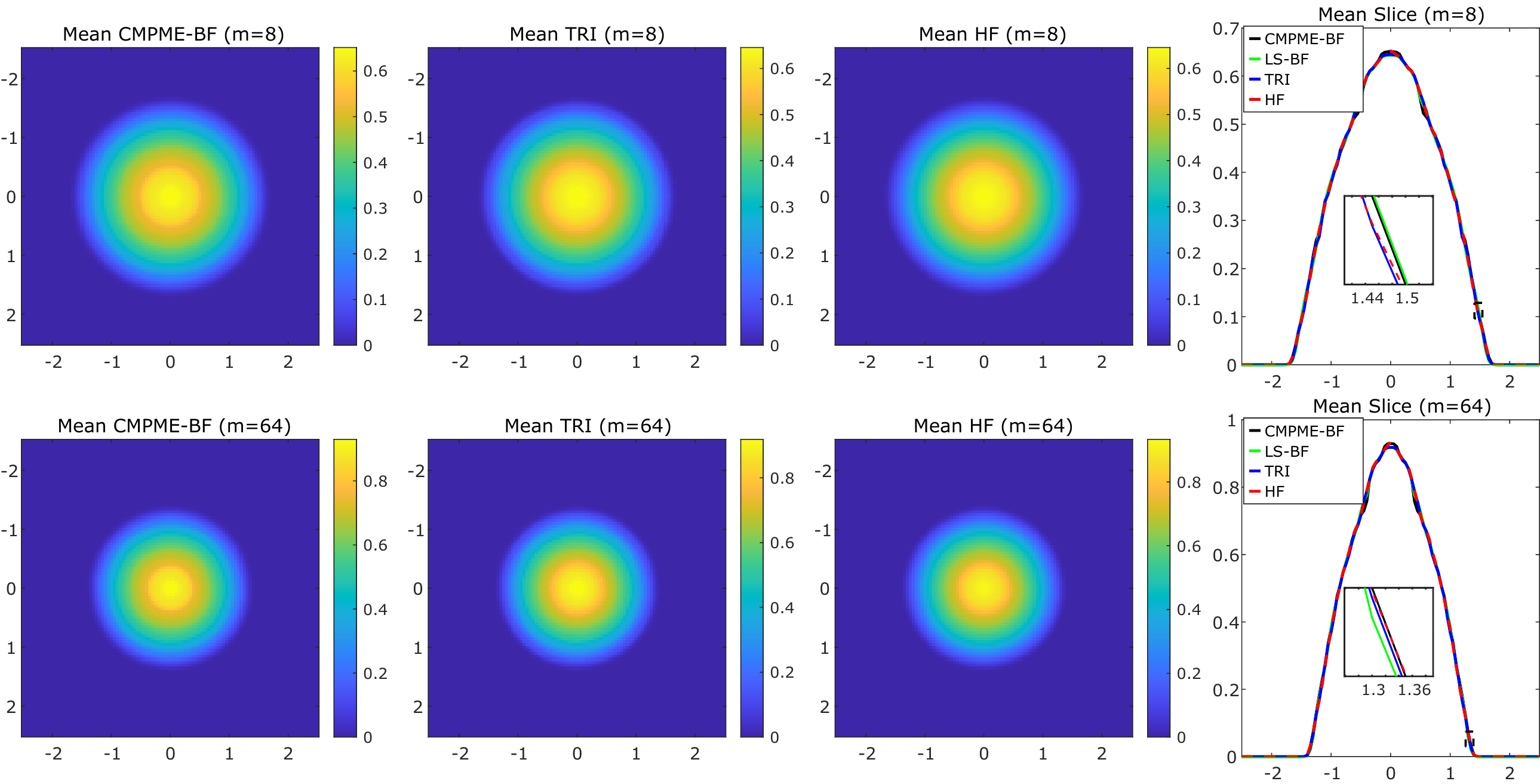}
    \caption{
Test 2(b): Mean field of \(\rho\) for CMPME-BF and TRI.
Top: \(m=8\); bottom: \(m=64\).
Columns show CMPME-BF, TRI reconstruction, HF reference, and
slices at $y=0$, which also include LS-BF.
}
    \label{fig:test2b_tri_mean}
\end{figure}

\begin{figure}[!htbp] 
    \centering
    \includegraphics[width=0.82\textwidth]{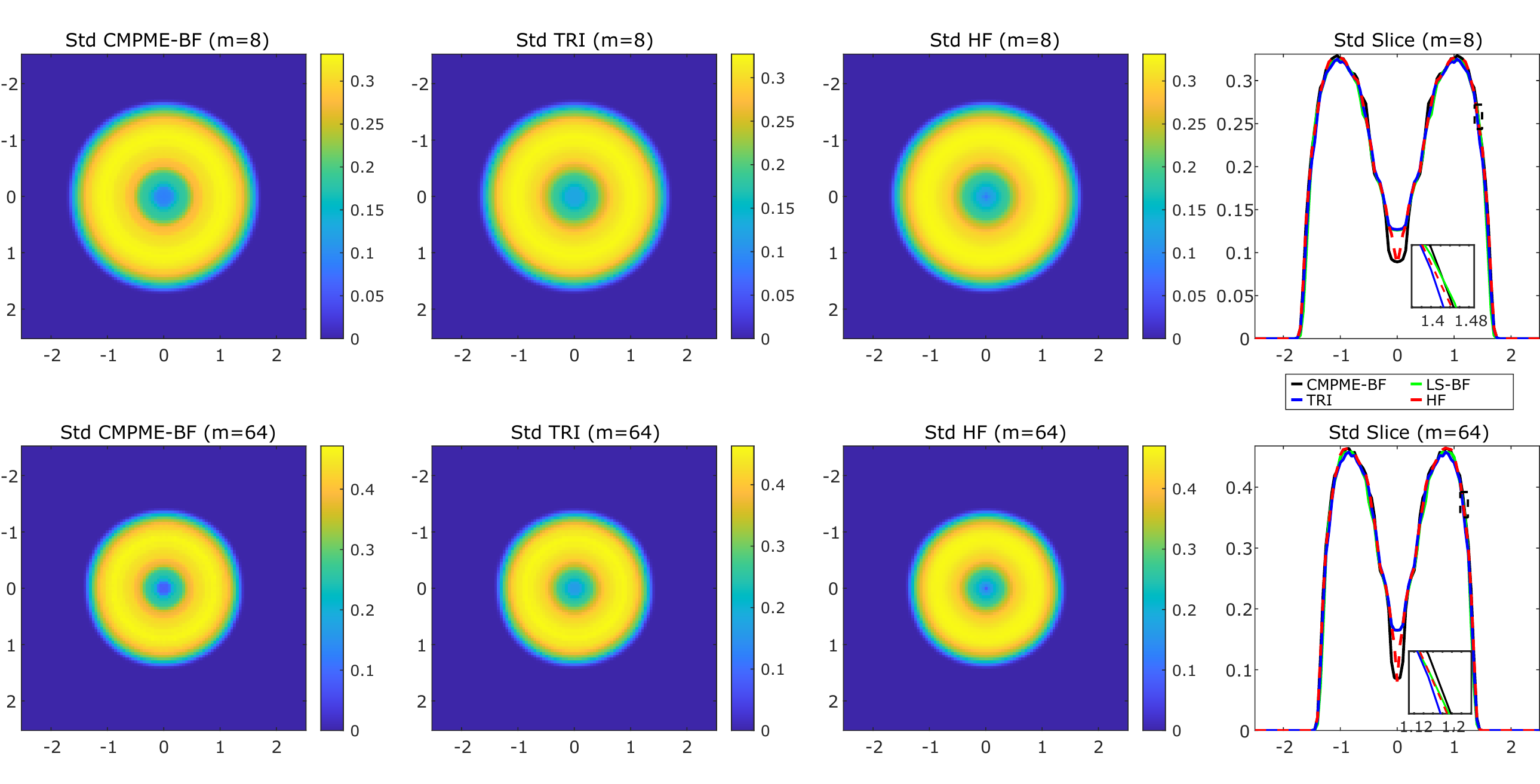}
    \caption{
Test 2(b): Standard deviation of \(\rho\) for CMPME-BF and TRI.
Top: \(m=8\); bottom: \(m=64\).
Columns show CMPME-BF, TRI reconstruction, HF reference, and
slices at $y=0$, which also include LS-BF.
}
    \label{fig:test2b_tri_std}
\end{figure}

\begin{figure}[!htbp] 
    \centering
    \includegraphics[width=0.82\textwidth]{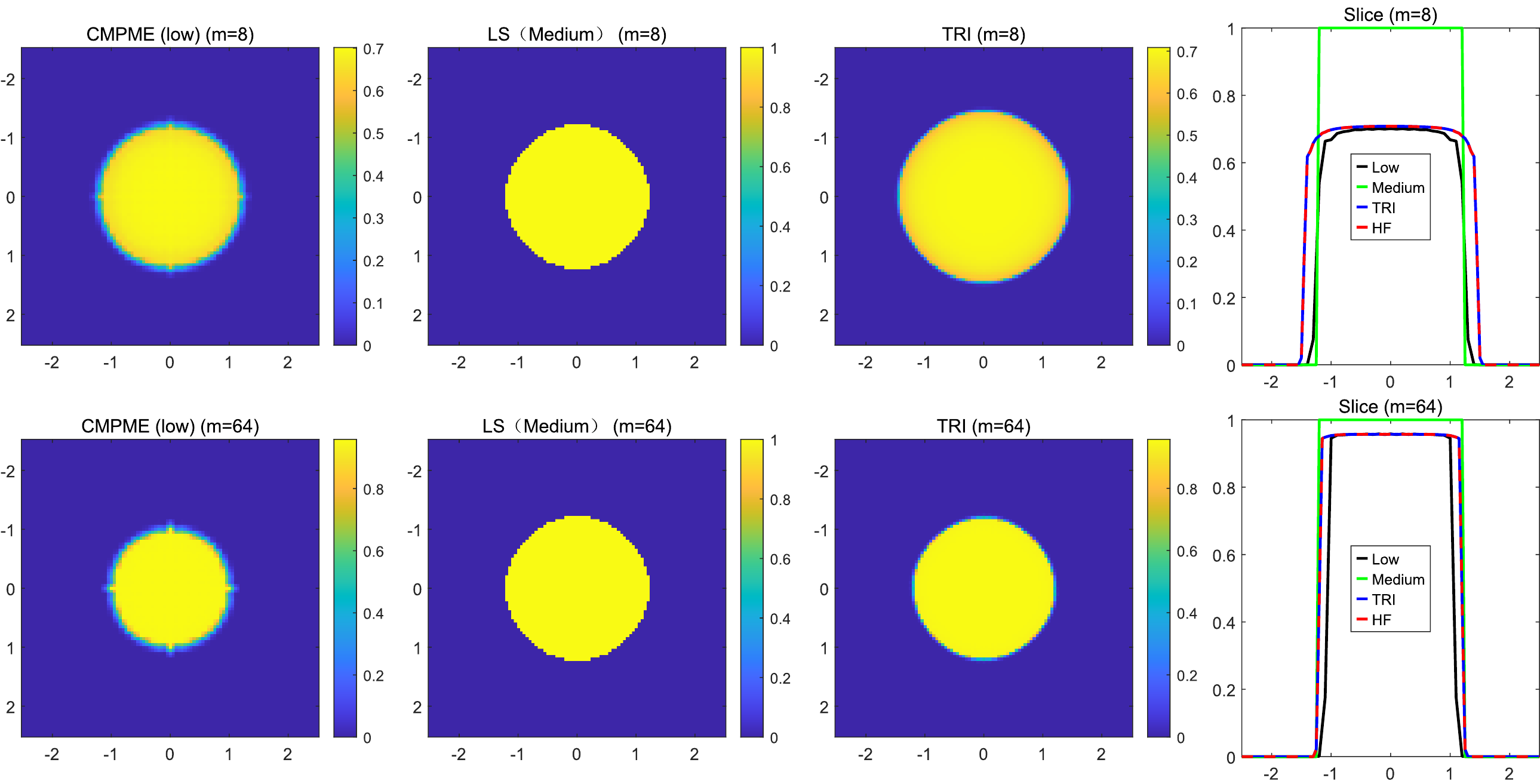}
    \caption{
Test 2(b): Representative realization of \(\rho\) at a fixed parameter $z$
under TRI approximation.
Top: \(m=8\); bottom: \(m=64\).
Columns show the CMPME low-fidelity solution, the LS intermediate-fidelity
solution, the TRI reconstruction, and one-dimensional slices at $y=0$,
including the HF reference.}
    \label{fig:test2ab_single}
\end{figure}

At $K=18$, TRI has the smallest reconstruction error for both
values of $m$ in Figure~\ref{fig:error_test2ab}. All three error estimators
remain above their corresponding observed errors and capture the error
scale. Since LS-BF statistics were presented in part~(a),
Figures~\ref{fig:test2b_tri_mean} and~\ref{fig:test2b_tri_std} focus on
CMPME-BF and TRI. The enlarged slices show how TRI improves the mean and
standard-deviation profiles near the interface while retaining their global
radial structure. For the fixed parameter in Figure~\ref{fig:test2ab_single},
CMPME represents a diffuse density profile, whereas LS represents the tumor
by a saturated region with a sharp boundary. The TRI reconstruction closely
recovers the HF interior density and interface transition for both values of $m$.

The reported LS time is approximately $19.686$ s per realization.
Using the $m=8$ CMPME and HF PME times reported in Test~1 as references
gives timing ratios of about $2.36$ for LS/CMPME and $16.17$ for HF/LS.
The improvement of TRI in Test~2(b) is obtained at the same HF basis
size ($K=18$) as the two BF methods.

\subsection{Test 3: Initial-interface complexity under oscillatory perturbations}

Following the radial perturbation setting in \cite{feng2023tumor},
we consider a lower-frequency mode ($n_{\rm petal}=6$) and a higher-frequency
mode ($n_{\rm petal}=16$) to examine frequency-dependent interface damping.
The radius and perturbation amplitude are random. LS supplies the LF model
for the baseline BF and RFPS-BF comparisons.
The model parameters are
\[
\lambda_i=\lambda_0(1+z_{1,i}),\quad
c_{B,i}=c_{B,0}(1+z_{1,i}),\quad
G_{0,i}=G_0,
\]
\[
R_{0,i}=R_0\Bigl(1+\sum_{k=1}^{d_z}\frac{z_{k,i}}{2k}\Bigr),
\qquad
A_i=A_0\Bigl(1+\sum_{k=1}^{d_z}\frac{z_{k,i}}{2k}\Bigr),
\]
where \(R_0=1\), \(A_0=0.1\), $\lambda_0=50$, $c_{B,0}=10$,
and $G_0=0.5$. The perturbed interface is defined as
\[
R_{\mathrm{petal}}(x,y;z_i)
=
R_{0,i}\bigl[1+A_i\cos(n_{\mathrm{petal}}\varphi)\bigr],
\]
where \(\varphi=\operatorname{atan2}(y,x)\) is the polar angle. The initial density is
\[
\rho_0(x,y;z_i)=
\begin{cases}
0.95, & \sqrt{x^2+y^2}<R_{\mathrm{petal}}(x,y;z_i),\\
0, & \text{otherwise}.
\end{cases}
\]

\begin{figure}[!htbp] 
    \centering
    \includegraphics[width=0.8\textwidth]{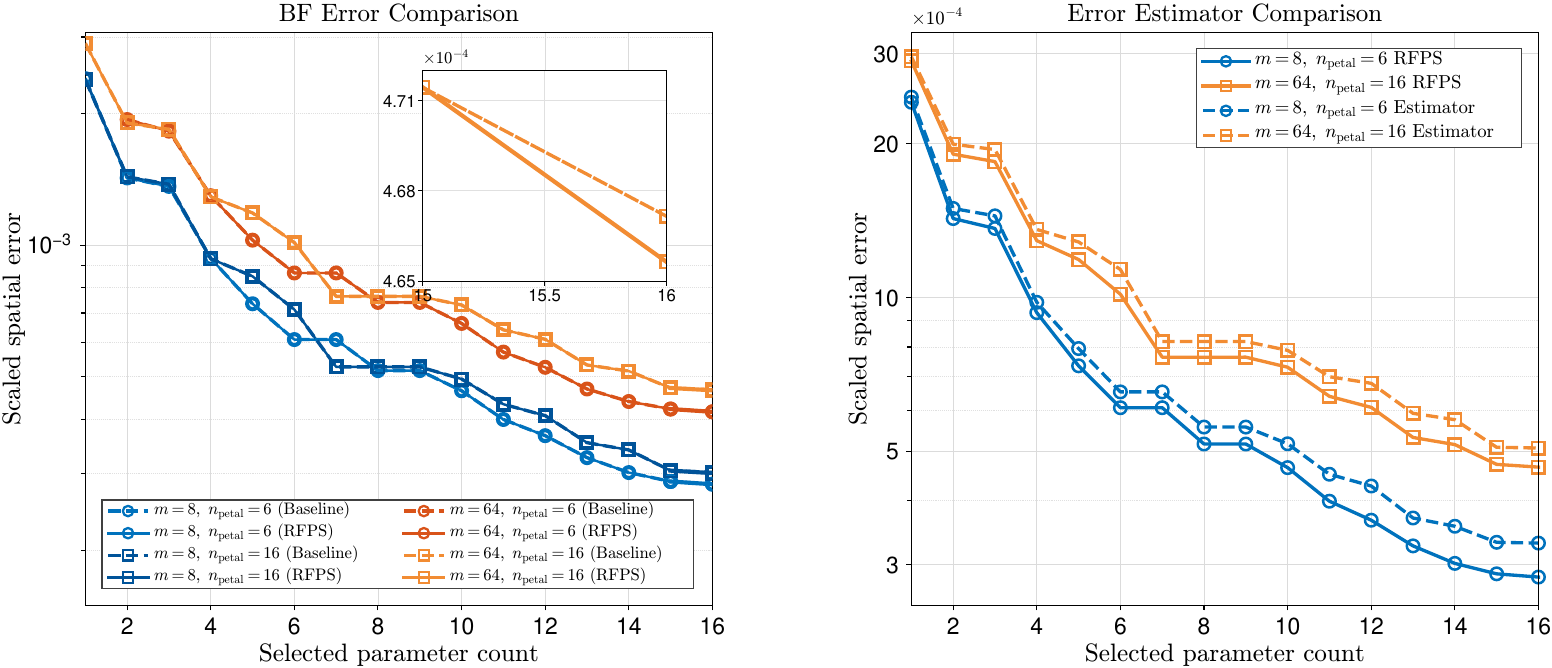}
    \caption{Test 3: Mean scaled spatial error on the parameter-matched
evaluation set versus the number of selected parameters.
Left: baseline BF and RFPS-BF for
$(m,n_{\rm petal})=(8,6),(8,16),(64,6),(64,16)$.
Right: RFPS-BF errors and error estimators for $(8,6)$ and $(64,16)$.}
    \label{fig:test3_error}
\end{figure}

For both petal counts and both values of $m$, RFPS reaches
$J\le0.5J_0$ after one addition, increasing the number of selected
parameters from 15 to 16. The baseline selects one additional sample by maximizing
the same inner-product score. In
Figure~\ref{fig:test3_error}, RFPS gives a modest further reduction for
$(m,n_{\rm petal})=(64,16)$ after the error curves nearly plateau.
The estimators for the two displayed cases remain above the observed
mean errors and track their scale.

\begin{figure}[!htbp]
    \centering
    \includegraphics[width=0.85\textwidth]{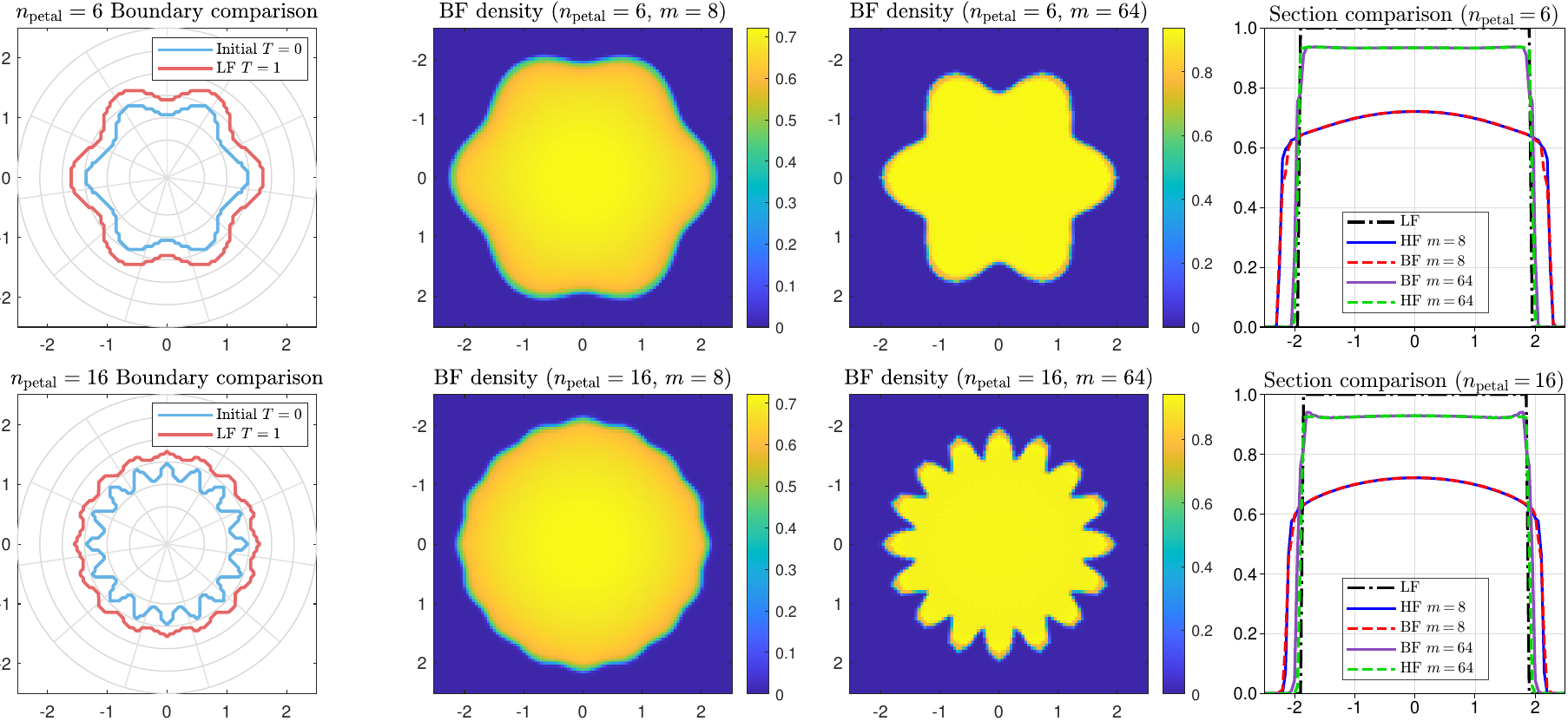}
    \caption{
Test 3: Representative interface evolution and reconstruction.
Top: \(n_{\mathrm{petal}}=6\); bottom: \(n_{\mathrm{petal}}=16\).
Left: illustrative LS interface growth from $t=0$ to $t=T=1$ for a fixed input.
The next two columns show initial-basis BF density fields at $m=8$ and
$m=64$ for an evaluation realization. Right: the corresponding BF and HF
slices at $y=0$, together with the LS profile.
}
    \label{fig:test3_single}
\end{figure}

The left panels of Figure~\ref{fig:test3_single} illustrate
frequency-dependent behavior consistent with \cite{feng2023tumor}:
as the tumor expands, the higher-frequency perturbation is damped more
strongly, while the six-lobed shape remains visible. The BF fields show
smoother interfaces for $m=8$ and more pronounced oscillations for $m=64$.
The corresponding slices closely match the HF interior density and edge
transitions.

\begin{figure}[!htbp] 
    \centering
    \includegraphics[width=0.8\textwidth]{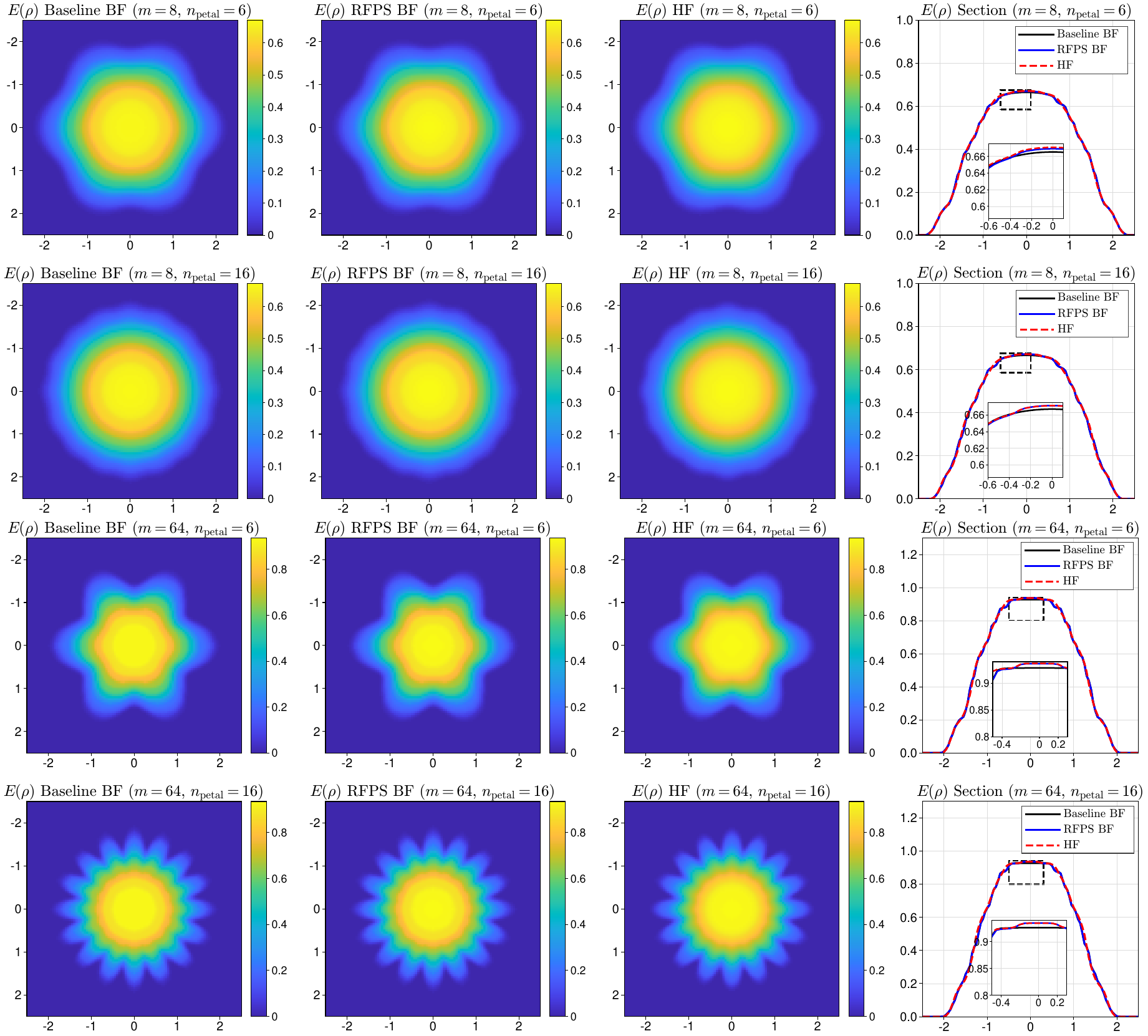}
    \caption{
Test 3: Mean field of \(\rho\) under oscillatory initial-interface
perturbations. Rows from top to bottom correspond to
\((m,n_{\mathrm{petal}})=(8,6),(8,16),(64,6),(64,16)\).
Columns show baseline BF, RFPS-BF, the HF reference, and slices
at $y=0$.
}
    \label{fig:test3_mean}
\end{figure}

\begin{figure}[!htbp]
    \centering
    \includegraphics[width=0.8\textwidth]{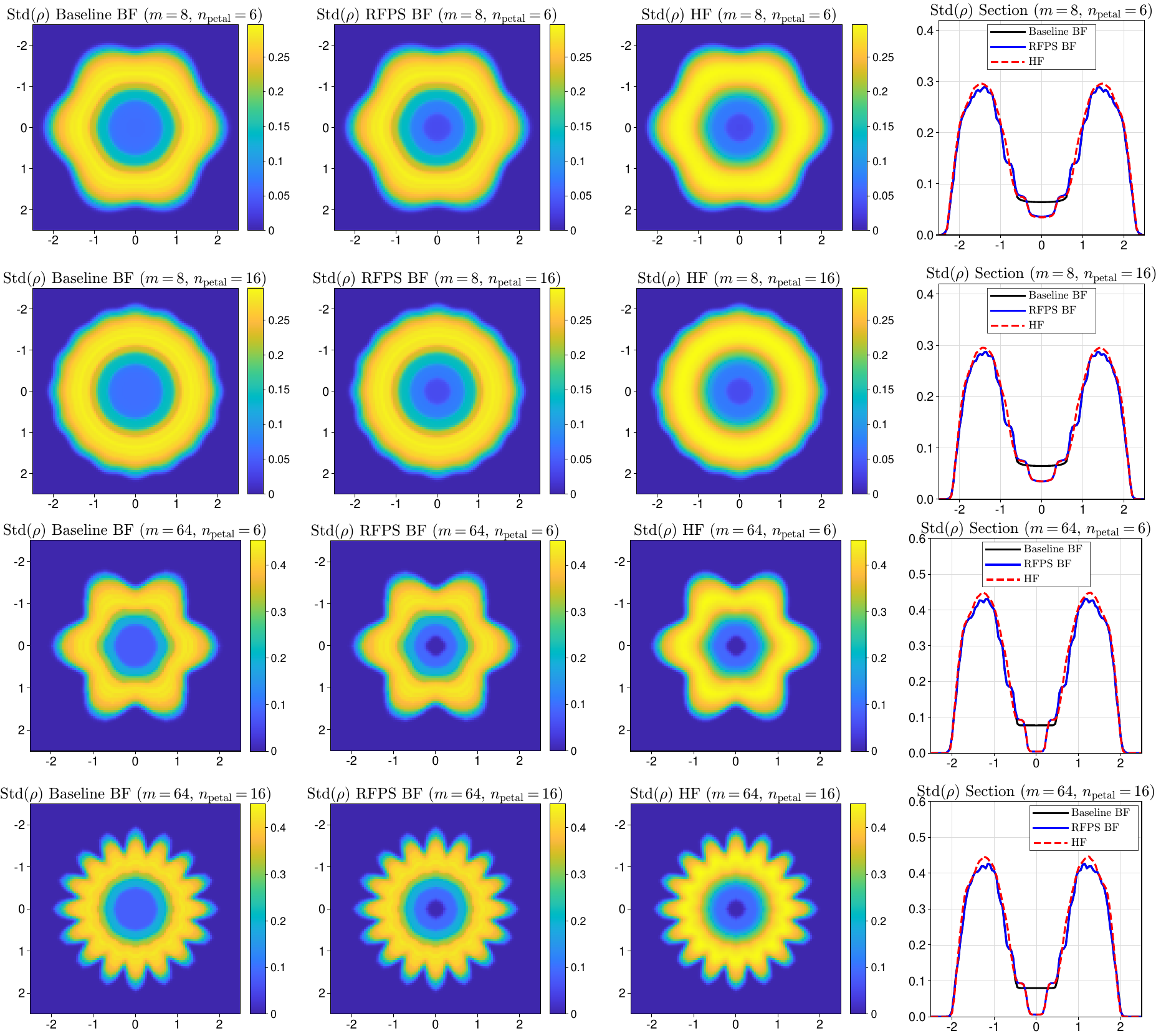}
    \caption{
Test 3: Standard deviation of \(\rho\) under oscillatory initial-interface
perturbations. Rows from top to bottom correspond to
\((m,n_{\mathrm{petal}})=(8,6),(8,16),(64,6),(64,16)\).
Columns show baseline BF, RFPS-BF, the HF reference, and slices
at $y=0$.
}
    \label{fig:test3_std}
\end{figure}

Figure~\ref{fig:test3_mean} shows that RFPS preserves the global mean
morphology and resolves the small variations in the central density profile,
as highlighted by the insets. In Figure~\ref{fig:test3_std}, RFPS retains
the peripheral uncertainty pattern and recovers the central depression
that is too shallow in the baseline, particularly for $m=64$. These
improvements are obtained at the same HF basis size as the baseline.

\FloatBarrier
\section{Conclusion}
We developed a hierarchical multi-fidelity method for uncertainty
quantification of porous-medium tumor growth models with moving interfaces.
The framework combines fine-grid asymptotic-preserving PME simulations
with coarse-grid PME and a level-set approximation of the Hele--Shaw limit,
exploiting both discretization-based and asymptotic-model-based fidelity
reduction.

The two-stage sampling procedure supplements pivoted Cholesky initialization
with RFPS, combining spatially weighted fitting residuals and farthest-point
separation to select additional HF sample locations. We established the
stability of the fixed-basis fit, quantified the decrease of its regularized
objective after enrichment, and related the residual stopping criterion
to HF reconstruction error under paired-model compatibility. An empirical
error estimator provides a common assessment for BF and TRI.

In the tested uncertainty settings, RFPS improves local density statistics,
particularly the standard-deviation profiles, after one to three
additions to the initial basis. The TRI construction achieves the smallest
final error in the nutrient-related hierarchy comparison at the same HF
basis size. The empirical estimators follow the observed reconstruction-error scales.
Future work will consider three-dimensional tumor growth \cite{liu2025three},
more complex microenvironmental coupling, and statistical error assessment.

\subsection*{Acknowledgments}

L. Liu acknowledges support from the National Key R\&D Program of China
(2021YFA1001200), the Ministry of Science and Technology of China
(co-PI, 2021), and the General Research Fund of the
Research Grants Council of Hong Kong (14303022, 2022; 14301423, 2023;
14307125, 2025; 14310726, 2026). Y. Feng was partially supported by the National Key
R\&D Program of China (2021YFA1001200) and the National Natural Science
Foundation of China (NSFC Youth Program, 12501669).

The HF Monte Carlo reference computations used resources of the
Supercomputing Center of Wuhan University. The authors also thank
Prof.~Ning Jiang (Wuhan University) and Dr.~Shuo Ling
(Shanghai Jiao Tong University) for valuable discussions and support.

\section*{Data and code availability}
The data and code in this article are available from the corresponding author upon reasonable request.

\section*{Conflict of interest}
The authors declare that they have no conflict of interest.

\FloatBarrier
\begin{appendices}

\section{Two-stage sampling details}
\label{app:sampling}
This appendix specifies the numerical settings and updates in
Algorithm~\ref{alg:two_stage_greedy}. The matrices $S_L$, $F_{\mathcal C}$
and weights $W$ are defined in Section~\ref{subsec:Two-stage greedy strategy}.

\textbf{Stage I: pivoted-Cholesky initialization.}
Initialize the remaining Gram diagonal by $q_i=f_i^Tf_i$ and the
permutation by $\sigma=(1,\dots,N)$. At step $j$, choose the largest
remaining $q_i$, swap it into position $j$, and apply the same swap
to the solution columns, $\sigma$, and the previously computed rows of
the Cholesky factor. The updates are
\[
 L_{jj}=\sqrt{q_j},\qquad
 L_{ij}=\frac{f_i^Tf_j-\sum_{\ell<j}L_{i\ell}L_{j\ell}}{L_{jj}},
 \qquad q_i\gets q_i-L_{ij}^2\quad(i>j).
\]
Here $f_i$ denotes the currently permuted column. Stop after
$\min\{K_0,N\}$ selections or when the largest remaining Gram diagonal
is below $10^{-8}$ (below $10^{-10}$ for TRI coarse-PME initialization).
Reset $K_0$ to the number retained. The original labels are
$\sigma(1),\dots,\sigma(K_0)$. For BF, the selected LF columns form
$\Phi_0$; for TRI, the same labels identify the LS columns forming
$\Phi_0$. Stage~II assumes $K_0\ge1$.

\textbf{Stage II: residual-based farthest-point sampling.}
Use the regularized fit and residual scores defined in
Section~\ref{subsec:Two-stage greedy strategy}, with $\delta=10^{-8}$.
The initial nearest-candidate representatives $j_a$ initialize the
minimum distances, not the selected parameter labels. With
$\epsilon=10^{-12}$, normalize the scores over the fixed candidate pool:
\[
 \widetilde r_i=\frac{r_i^{\rm sel}-\min_jr_j^{\rm sel}}
 {\max\{\max_jr_j^{\rm sel}-\min_jr_j^{\rm sel},\epsilon\}},
 \qquad
 \widetilde d_i=\frac{d_i}{\max\{\max_jd_j,\epsilon\}}.
\]
The RFPS ranking is $\omega\widetilde r_i+(1-\omega)\widetilde d_i$,
maximized over the unexamined indices. Ties in pool searches are resolved
by pool order. Before appending a candidate $i$, require
\[
 r_i^{\rm sel}\ge10^{-3},\qquad
 \|\mathcal R(:,i)\|_2\ge10^{-8},\qquad
 \max_a|\operatorname{corr}(g_i,\Phi(:,a))|\le0.999,
\]
where $\operatorname{corr}$ is the Pearson correlation across spatial
entries, defined for columns of nonzero spatial variance.
An examined candidate, including a rejected one, is not reconsidered.
Rejection leaves the basis and distances unchanged. Acceptance appends
$g_i$, records its parameter label for the HF solve, and updates the
stored distances of unexamined candidates by the minimum-distance rule
in Algorithm~\ref{alg:two_stage_greedy}. The distances of examined
candidates are left unchanged for the normalization. The budget
$K_{\rm add}$ counts accepted additions rather than examined candidates.

For the stopping test, compute
$J(\Phi)=\|W^{1/2}\mathcal R\|_F/\sqrt{N_CS_w}$ over the entire fixed
pool, including previously examined candidates. Check $J(\Phi)\le\theta J_0$
before each selection and after the last accepted addition. If the target
is unmet, report whether the addition budget or pool was exhausted.

\section{Asymptotic-preserving PME discretization}
\label{app:PME}

This appendix summarizes the main discretization steps of the HF solver
based on the prediction--correction asymptotic-preserving framework
of Section~\ref{sec:solver}. At each time level, $c^n$ is obtained
from the nutrient problem \eqref{eq:c_vivo} with the current density
$\rho^n$; the velocity predictor and density update then use this field.
\vspace{2mm}

\textbf{Implicit--explicit discretization. }
All variables are stored at uniform grid points
\((x_i,y_j)\in[a,b]^2\), with spatial step sizes \(\Delta x,\Delta y\).
The components $u$ and $v$ denote horizontal and vertical
velocities. The half-grid densities in the velocity predictor use
arithmetic averages.
The implicit--explicit (IMEX) scheme treats the predicted velocity
implicitly, with density-dependent coefficients frozen at the previous
time level. The density update uses explicit fluxes and a semi-implicit
reaction term.
\vspace{2mm}

\textbf{Prediction step: velocity update. }
The IMEX update for the horizontal velocity component \(u^{n*}\) reads
\begingroup\small
\begin{equation}\label{eq:HF_u_update_app}
\begin{aligned}
\frac{u^{n*}_{ij}-u^n_{ij}}{\Delta t}
&=
m\Bigg[
\frac{
(\rho^n_{i+1/2,j})^{m-2}
\big((\rho^nu^{n*})_{i+1,j}-(\rho^nu^{n*})_{i,j}\big)
-
(\rho^n_{i-1/2,j})^{m-2}
\big((\rho^nu^{n*})_{i,j}-(\rho^nu^{n*})_{i-1,j}\big)
}{\Delta x^2}
\\[-2pt]
&\quad+
\frac{
(\rho^n_{i+1,j})^{m-2}
\big((\rho^nv^{n*})_{i+1,j+1}-(\rho^nv^{n*})_{i+1,j-1}\big)
-
(\rho^n_{i-1,j})^{m-2}
\big((\rho^nv^{n*})_{i-1,j+1}-(\rho^nv^{n*})_{i-1,j-1}\big)
}
{4\Delta x\,\Delta y}
\\[-2pt]
&\quad-
\frac{
(\rho^n_{i+1,j})^{m-2}G_0c^n_{i+1,j}\rho^n_{i+1,j}
-
(\rho^n_{i-1,j})^{m-2}G_0c^n_{i-1,j}\rho^n_{i-1,j}
}{2\Delta x}
\Bigg].
\end{aligned}
\end{equation}
\endgroup
The update for the vertical component \(v^{n*}\) is obtained analogously
by interchanging $(x,i,u)$ and $(y,j,v)$.
\vspace{2mm}

\textbf{Density update. } The face states $\rho_L$ and $\rho_R$
are reconstructed from $\rho^n$ with the minmod limiter; $u^*$ and $v^*$
are arithmetic face averages of the predicted nodal velocities
$u^{n*}$ and $v^{n*}$.
Using central-upwind numerical fluxes,
\[
F_{1,i\pm1/2,j}
=
\frac12
\left(
\rho_L u^{*}
+
\rho_R u^{*}
-
|u^{*}|(\rho_R-\rho_L)
\right),
\]
and
\[
F_{2,i,j\pm1/2}
=
\frac12
\left(
\rho_L v^{*}
+
\rho_R v^{*}
-
|v^{*}|(\rho_R-\rho_L)
\right),
\]
the density is updated by
\begin{equation}\label{eq:HF_density_app}
\frac{\rho^{n+1}_{ij}-\rho^n_{ij}}{\Delta t}
+
\frac{F_{1,i+1/2,j}-F_{1,i-1/2,j}}{\Delta x}
+
\frac{F_{2,i,j+1/2}-F_{2,i,j-1/2}}{\Delta y}
=
G_0 c^n_{ij}(z)\rho^{n+1}_{ij}.
\end{equation}
The flux update is conservative, and the reaction term changes the
total mass. The semi-implicit update requires $1-\Delta t G_0c^n_{ij}>0$
for a positive reaction denominator.
\vspace{2mm}

\textbf{Correction step. }
The correction step enforces the porous-medium constitutive relation by
resetting the velocity field to the corresponding discrete gradient
form:
\[
\begin{aligned}
u^{n+1}_{ij}
&=
-\frac{m}{m-1}
\frac{
(\rho^{n+1}_{i+1,j})^{m-1}
-
(\rho^{n+1}_{i-1,j})^{m-1}
}{2\Delta x},
\\[2pt]
v^{n+1}_{ij}
&=
-\frac{m}{m-1}
\frac{
(\rho^{\,n+1}_{i,j+1})^{m-1}
-
(\rho^{\,n+1}_{i,j-1})^{m-1}
}{2\Delta y}.
\end{aligned}
\]
This step restores consistency with the nonlinear diffusion constraint
and completes one HF time update.

\section{Proofs of stability and accepted enrichment}
\label{app:analysis_proofs}

\subsection{Fixed-basis stability}
\label{app:ridge_proof}
\begin{proof}[Proof of Lemma~\ref{lem:ridge_stability}]
For every nonzero vector $v$,
$v^T(\Phi^T\Phi+\delta I)v=
\|\Phi v\|_2^2+\delta\|v\|_2^2>0$.
Thus the objective
$\|F_{\mathcal C}-\Phi C\|_F^2+\delta\|C\|_F^2$ is strictly convex,
and its normal equations have the unique solution
$C_\delta=A_\delta F_{\mathcal C}$.
For a thin singular value decomposition $\Phi=U\Sigma V^T$,
$A_\delta=V\operatorname{diag}
\{\sigma_j/(\sigma_j^2+\delta)\}U^T$, including zero singular values
when the basis is rank deficient.
The inequality $\sigma^2+\delta\ge2\sigma\sqrt\delta$ proves
\eqref{eq:ridge_stability}. Since
$\Delta C_\delta=A_\delta\Delta F$, the matrix norm inequality
$\|AB\|_F\le\|A\|_2\|B\|_F$ proves the perturbation estimate;
applying it once more to $\Psi\Delta C_\delta$ proves the HF statement.
\end{proof}

\subsection{Accepted enrichment}
\label{app:descent_proof}
\begin{proof}[Proof of Proposition~\ref{prop:accepted_descent}]
At the minimizer $C_\delta$, the normal equations give
$\mathcal V_\delta(\Phi)=\|F_{\mathcal C}\|_F^2-
\operatorname{tr}(F_{\mathcal C}^T\Phi A_\delta F_{\mathcal C})$.
By the definition of $Q_\delta$, $\mathcal V_\delta(\Phi)=
\operatorname{tr}(F_{\mathcal C}^TQ_\delta F_{\mathcal C})$.
Multiplying $Q_\delta$ by $\Phi\Phi^T+\delta I$ gives $\delta I$, so
\begin{equation}
\label{eq:q_delta_identity}
Q_\delta=\delta(\Phi\Phi^T+\delta I)^{-1}.
\end{equation}
Set $M=\Phi\Phi^T+\delta I\succ0$. Applying the Sherman--Morrison
formula gives
$(M+aa^T)^{-1}=M^{-1}-M^{-1}aa^TM^{-1}/(1+a^TM^{-1}a)$.
Multiplying by $\delta$ and using \eqref{eq:q_delta_identity} yields
\[
Q_{\delta,+}=Q_\delta-
\frac{Q_\delta aa^TQ_\delta}{\delta+a^TQ_\delta a}.
\]
Subtracting the two trace expressions gives the equality in
\eqref{eq:objective_decrease}.

For the quantitative bound, \eqref{eq:q_delta_identity} gives
$0\prec Q_\delta\preceq I$ and $Q_\delta^2\preceq Q_\delta$.
Since $r_i^{\rm sel}=\|W^{1/2}Q_\delta a\|_2/\sqrt{S_w}$,
putting $t=a^TQ_\delta a$ yields
\[
t\ge\|Q_\delta a\|_2^2
\ge\frac{\|W^{1/2}Q_\delta a\|_2^2}{w_{\max}}
=\alpha_a>0.
\]
Because $a$ is a column of $F_{\mathcal C}$,
$\|F_{\mathcal C}^TQ_\delta a\|_2^2\ge t^2$.
The function $t\mapsto t^2/(\delta+t)$ is increasing for $t\ge0$,
which proves the lower bound in \eqref{eq:objective_decrease}.
The acceptance threshold in Appendix~\ref{app:sampling} additionally
gives $\alpha_a\ge S_w10^{-6}/w_{\max}$.
Finally,
$\|Q_\delta F_{\mathcal C}\|_F^2
\le\operatorname{tr}(F_{\mathcal C}^TQ_\delta F_{\mathcal C})$
by $Q_\delta^2\preceq Q_\delta$, and
$\|W^{1/2}Q_\delta F_{\mathcal C}\|_F^2
\le w_{\max}\|Q_\delta F_{\mathcal C}\|_F^2$.
Dividing by $N_CS_w$ proves the stated envelope for $J^2$.
\end{proof}

\end{appendices}

\bibliographystyle{plain}
\bibliography{main}
\end{document}